\def\no{\if01}
\def\iftwelvept{\no}

\def\ifusepdf{\no}
\def\ifpsfont{\no}

\iftwelvept
\documentclass[leqno,12pt]{amsart}
\else
\documentclass[leqno,12pt]{amsart}
\fi
\usepackage{amssymb}
\usepackage{amscd}
\usepackage{latexsym}
\usepackage{verbatim}
\usepackage{mathrsfs}
\usepackage[all]{xy}

\ifusepdf
\usepackage{hyperref}
\else\fi
\ifpsfont
\usepackage[T1]{fontenc}
\usepackage{times}
\else\fi

\iftwelvept
\else\fi

\theoremstyle{plain}
\newtheorem{Theorem}{Theorem}[section]

\newtheorem{Proposition}[Theorem]{Proposition}
\newtheorem{Lemma}[Theorem]{Lemma}
\newtheorem{Corollary}[Theorem]{Corollary}

\newtheorem{Claim}{Claim}[Theorem]

\theoremstyle{definition}

\newtheorem{Definition}[Theorem]{Definition}
\newtheorem{Remark}[Theorem]{Remark}
\newtheorem{Example}[Theorem]{Example}

\renewcommand{\theTheorem}{\arabic{section}.\arabic{Theorem}}
\renewcommand{\theClaim}{\arabic{section}.\arabic{Theorem}.\arabic{Claim}}
\renewcommand{\theequation}{\arabic{section}.\arabic{Theorem}.\arabic{Claim}}

\newcommand{\ZZ}{{\mathbb{Z}}}
\newcommand{\QQ}{{\mathbb{Q}}}

\newcommand{\CC}{{\mathbb{C}}}

\newcommand{\GG}{{\mathbb{G}}}

\newcommand{\YY}{{\mathcal{Y}}}

\newcommand{\II}{{\mathcal{I}}}
\newcommand{\FF}{{\mathcal{F}}}
\newcommand{\WW}{{\mathcal{W}}}

\newcommand{\OO}{{\mathcal{O}}}

\newcommand{\LL}{{\mathsf{L}}}

\newcommand{\GGG}{\mathcal{G}}

\newcommand{\UU}{\mathcal{U}}

\newcommand{\Hom}{\operatorname{Hom}}
\newcommand{\Ext}{\operatorname{Ext}}

\newcommand{\XPre}{X(\textup{Pre})}

\newcommand{\Stab}{\operatorname{Stab}}

\newcommand{\Spec}{\operatorname{Spec}}
\newcommand{\Spf}{\operatorname{Spf}}

\newcommand{\Isom}{\underline{\operatorname{Isom}}}

\newcommand{\gs}{\operatorname{gs}}

\newcommand{\Aut}{\underline{\operatorname{Aut}}}

\newcommand{\XX}{\mathcal{X}}

\newcommand{\Proof}{{\sl Proof.}\quad}
\newcommand{\QED}{{\unskip\nobreak\hfil\penalty50\quad\null\nobreak\hfil
{$\Box$}\parfillskip0pt\finalhyphendemerits0\par\medskip}}

\begin{document}

\title{Stable points on algebraic stacks}
\thanks{E-mail address: iwanari@kurims.kyoto-u.ac.jp (Corresponding author) \\\ \ \ \ The author is supported by Grand-in-Aid for JSPS}
\author[Isamu Iwanari]{Isamu Iwanari \\\    \\ Research Institute for Mathematical Sciences, Kyoto University, Kyoto, 606-8502, Japan }

\begin{abstract}
This paper is largely concerned with constructing coarse moduli spaces
for Artin stacks.
The main purpose of this paper is to introduce the notion of stability on
an arbitrary Artin stack and construct a coarse moduli space for the open substack
of stable points.
Also, we present an application to coherent cohomology of Artin stacks.
\vspace{1mm}

\begin{flushleft}
\end{flushleft}

\end{abstract}
\maketitle


\section*{Introduction}
\renewcommand{\theTheorem}{\Alph{Theorem}}

This paper is largely concerned with constructing coarse moduli spaces
for Artin stacks.
Roughly speaking, a coarse moduli space for an Artin stack $\XX$ is
the best approximation of $\XX$ by an algebraic space
whose underlying space coincides with that of $\XX$ (see section 1).
Coarse moduli spaces play the role of the bridge between the geometry
of Artin stacks and that of algebraic spaces and schemes.
Even if we are ultimately interested in stacks,
the existence of coarse moduli spaces is useful in various situations.
In their influential paper (\cite{KM}),
Keel and Mori proved the existence of a coarse moduli space for
an Artin stack whose objects have finite automorphism groups.
In particular, their theorem implies the existence of a coarse
moduli space for a Deligne-Mumford stack under a weak assumption.
Let us turn our attention to arbitrary Artin stacks.
There are many Artin
stacks which have positive dimensional automorphisms.
For instance, such examples arise from group actions on algebraic spaces,
moduli spaces of vector bundles and complexes on algebraic varieties, affine geometry and so on.
Hence it is desired to construct coarse moduli spaces for general
Artin stacks.
However, we can readily find Artin stacks
which do not admit coarse moduli spaces.
For example, the quotient stack $[\mathbb{A}^1/\GG_m]$ arising from the natural
action
of the torus $\GG_m\subset \mathbb{A}^1$ on the affine line
$\mathbb{A}^1$ does not have a coarse moduli space.
Thus to construct coarse moduli spaces, we need to impose 
some condition on Artin stacks.

Let $X$ be an algebraic scheme and $G$ a reductive group acting on $X$.
In his theory of Geometric Invariant Theory (\cite{GIT}), Mumford defined the notion of (pre-)stable points on $X$
and proved that the quotient of the open subset of (pre-)stable points
by $G$ exists as a geometric quotient.
In terms of stacks, it says that
the open substack of the quotient stack $[X/G]$ associated to the
open set of (pre-)stable
points has a ``coarse moduli scheme''.

 The main purpose of this paper is to introduce the notion of stability on
an arbitrary Artin stack and to construct a coarse moduli space for the open substack
of stable points.
By the universality of coarse moduli space
its existence depends only upon the local properties of each point.
In addition, one of the reasons why the theorem of Keel and Mori
is useful, is that the finiteness of automorphisms
is a local condition, and so it can be checked pointwisely.
Likewise, if $p$ is a point on an Artin stack $\XX$, then
the stability (of the point $p$) introduced in this paper is
defined by using local data around $p$.
The central notion introduced in this paper is GIT-like p-stability (cf. Definition~\ref{defgit}).
(The letter ``p'' stands for ``pointwise'' and ``potential''.)
In its naive form, the first main result of this paper states the
following:

\begin{Theorem}
\label{intromain}
Let $\XX$ be an Artin stack locally of finite type over a perfect field.
Then GIT-like p-stable points form an open substack
$\XX^{\gs}$ and there exists a coarse moduli map
\[
\pi: \XX^{\gs}\longrightarrow X
\]
such that $X$ is an algebraic space locally of finite type.
If $X'\to X$ is a flat morphism of algebraic spaces, then
the second projection $\XX^{\gs}\times_{X}X'\to X'$ is also
a coarse moduli map. $($See Theorem~\ref{maingit} for detail.$)$
Moreover, the followings hold:
for any point $x\in X$, there exists an \'etale neighborhood
$U\to X$ such that $\XX^{\gs}\times_{X}U$ has the form $[W/G]$ of a quotient
stack,
where $W$ is affine over $U$ and $G$ is a linearly reductive group.
$($cf. Corollary~\ref{gitlocal}.$)$

\end{Theorem}

Let us recall that in the case of Deligne-Mumford stacks the theorem of Keel and Mori implies a more precise and fruitful correspondence: a Deligne-Mumford stack $\XX$
has finite
inertia stack, i.e. every object in $\XX$ has a finite automorphism group scheme if and only if the following conditions hold:
(i) $\XX$ has a coarse moduli space $X$ and the formation of coarse moduli space
commutes with flat base change $X'\to X$, and (ii) for any point on $X$
there exists an \'etale neighborhood $U\to X$ such that $\XX\times_{X}U$
has the form $[V/G]$ of a quotient stack, where $V$ is affine over $U$ and
$G$ is a finite group acting on $V$ over $U$.
In characteristic zero, our notion of GIT-like p-stability
successfully generalizes this correspondence to Artin stacks:

\begin{Theorem}
\label{B}
Let $\XX$ be an Artin stack locally of finite type over a
field of characteristic zero. Let $\XX^{gs}$ be the open substack of
GIT-like p-stable points on $\XX$.
Then $\XX^{\gs}$ coincides with $\XX$ 
if and only if the following conditions hold:
\begin{enumerate}
\renewcommand{\labelenumi}{(\roman{enumi})}

\item $\XX$ has a coarse moduli space $X$ and the formation of coarse moduli space
commutes with flat base change $X'\to X$,

\item for any point on $X$ there exists an \'etale neighborhood
$U\to X$ such that $\XX\times_{X}U$ has the form $[V/G]$ of a quotient stack,
where $V$ is affine over $U$ and $G$ is a reductive group acting on $V$.
\end{enumerate}
$($See Theorem~\ref{characterization} for the precise statement.$)$
\end{Theorem}

It is worth mentioning that in characteristic zero
the definition of GIT-like p-stability is a natural generalization
of the finiteness of automorphism groups, and it is described
by the reductive condition on the automorphism group (cf. (a) in Definition~\ref{defgit})
and the conditions (a), (b), (c) in Remark~\ref{gitdefremark} (iii) (see also Remark~\ref{zerodef}).
In the Deligne-Mumford case, the existence of coarse moduli space
and the \'etale local quotient structure have played an important
role. For example, To\"en's Riemann-Roch theorem for Deligne-Mumford stacks
(\cite{Toen}) relies on them, and Gromov-Witten theory of Deligne-Mumford
stacks (\cite{AGV}) requires them.
Similarly, in our general situation coarse moduli spaces
and the \'etale local quotient structures are quite useful.
Indeed we present applications of Theorem A and B
to coherent cohomology (see section 7).

As an example of our situation, our existence theorem
contains the case of Geometric Invariant Theory.
Namely, the relationship with Geometric Invariant Theory is described
as follows:

\begin{Theorem}
\label{C}
Let $X$ be a scheme of locally
finite type and separated over an algebraically closed field $k$
of characteristic zero.
Let $G$ be a reductive group acting on $X$.
Let $\XPre$ be the open subset of $X$ consisting of pre-stable points
in the sense of Geometric Invariant Theory $($\cite{GIT}$)$.
Let $\mathcal{S}$ be the maximal open substack of the quotient stack
$[X/G]$ such that $\mathcal{S}$ admits a coarse moduli space which is
a scheme. Then we have
\[
[\XPre/G]=[X/G]^{\gs}\cap \mathcal{S}.
\]
$($See Theorem~\ref{gitrelation}, Remark~\ref{Crem}.$)$
\end{Theorem}


\vspace{2mm}

This paper is organized as follows.
In section 1, we recall some basic facts and fix some notation.
In section 2 and 3 we present preliminary notions and preparatory results.
We define strong p-stability and
 prove the existence of coarse moduli spaces in
 the case of strongly p-stable case.
In section 4, we introduce the method of deformations of
coarse moduli spaces, which we apply in section 5.
In section 5, we then introduce the notion of GIT-like p-stability,
which is the main notion of this paper.
By studying the local structure of a GIT-like p-stable point,
we prove Theorem A.
In section 6, we discuss the relationship with Geometric Invariant Theory
due to Mumford (in characteristic zero). Namely, we prove Theorem C.
From Theorem A and the comparison to Geometric Invariant Theory
we deduce Theorem B.
Finally, in section 7 we present applications
concerning coherent cohomology.
In Appendix, for the reader's convenience we collect some results
on limit arguments and rigidifications.

\vspace{3mm}

{\it Acknowledgements}.
I would like to thank H. Minamoto who helped me to prove some
lemmas and S. Mori who asked me the properties of coarse moduli maps
and it led me to section 7.
I would like to express my graditude to J. Alper
for his helpful comments and to the referee for valuable comments
and suggestions.
I would also like to acknowledge a great intellectual debt to the works of various authors, though it is difficult to fully acknowledge them here.
Finally, I am supported by JSPS grant.

\renewcommand{\theTheorem}{\arabic{section}.\arabic{Theorem}}
\renewcommand{\thesubsubsection}{\arabic{section}.\arabic{subsection}.\arabic{subsubsection}}

\section{Preliminaries}

We establish some notation, terminology and basic facts, which we will use in this paper.

\vspace{2mm}

We refer to the book \cite{LM} as the general reference to the notion
of algebraic stacks. We use henceforth
the word ``algebraic stack'' instead of Artin stack.
Similarly to this reference, in this paper except Appendix,
all schemes, algebraic spaces, and algebraic stacks
are assumed to be quasi-separated.
Let $S$ be a scheme.
For an algebraic stack $p:\XX\to S$,
the fiber $\XX(U)$ over an $S$-scheme $U$
is the subcategory whose objects is objects that lies over $U$
and whose morphisms are morphisms $f$ where $p(f)$ is the identity of $U$.
If $U=\Spec A$ we often write $\XX(A)$ for $\XX(U)$.
For any $S$-scheme $T$ and any $\xi, \eta:X\to \XX$,
let
\[
\Isom_{\XX,T}(\xi,\eta):(X\textup{-schemes})\to (\textup{sets})
\]
be the functor which to any $f:V\to T$ associates
the set of (iso)morphisms from $f^*\xi$ to $f^*\eta$ in $\XX(V)$.
Here we denote by $f^*\xi\in \XX(V)$ (resp. $f^*\eta$)
the pullback of $\xi$ (resp. $\eta$),
that is, the image $\xi(V\to X)\in \XX(V)$ (resp. $\eta(V\to T$)).
The functor $\Isom_{\XX,T}(\xi,\eta)$ is naturally represented by the fiber product
$X\times_{(\xi,\eta),\XX\times_S\XX,\Delta}\XX$,
where $\Delta$ is diagonal. In particular, it is an algebraic space separated
and of finite type over $T$.
If $\xi=\eta$, then $\Isom_{\XX,T}(\xi,\eta)$ is a group algebraic space
over $X$.
In this case, we write $\Aut_{\XX,T}(\xi)$ for
$\Isom_{\XX}(\xi,\xi)$.
If no confusion seems likely to arises, we abuse notation
and write $\Isom(\xi,\eta)$ or $\Isom_{\XX}(\xi,\eta)$ or $\Isom_X(\xi,\eta)$ for 
$\Isom_{\XX,T}(\xi,\eta)$.
If $G$ is a flat group scheme of finite and separated
over a scheme $S$,
then $BG$ denotes the classifying stack of principal $G$-bundles,
that is, the quotient stack $[S/G]$ associated to
the trivial action.

\vspace{2mm}

Let $\XX$ be an algebraic stack over a scheme $S$.
A {\it coarse moduli map} for $\XX$ over $S$ is
a morphism $\pi:\XX\to X$ from $\XX$ to an algebraic space $X$ over $S$
such that the following conditions hold.
\begin{enumerate}
\item If $K$ is an algebraically closed $S$-field,
then the map $\pi$ induces a bijective map between
the set of isomorphism classes of objects in $\XX(K)$ and $X(K)$.

\item The map $\pi$ is universal for maps from $\XX$ to algebraic spaces
over $S$.

\end{enumerate}
We shall refer to $X$ as a {\it coarse moduli space} for $X$.

Let us recall the theorem of Keel and Mori (\cite{KM}): Let $\XX$ be an
algebraic stack $\XX$ locally of finite type over a locally noetherian scheme $S$.
Assume that for any object $\xi\in \XX(T)$ over an affine $S$-scheme $T$
the automorphism group $\Aut_{\XX,T}(\xi)$ is finite over $T$.
Then there exists a coarse moduli map
\[
\pi:\XX\to X,
\]
where $X$ is an algebraic space locally of finite type over $S$.
Moreover, the natural homomorphism $\OO_{X}\to \pi_*\OO_{\XX}$
is an isomorphism, and
$\pi$ is quasi-finite and proper.
If $\XX$ has finite diagonal, then $X$ is separated over $S$.
We shall henceforth refer to this theorem as {\it Keel-Mori theorem}.

\section{Strong p-stability}

First of all, recall the notion of effective versal deformations (see \cite[section 1]{A1}, \cite[(10.11)]{LM}).
Let $\XX$ be an algebraic stack over a scheme $S$.
Let $K$ be an $S$-field, and let $\xi_0\in \FF(K)$.
An effective versal deformation of $\xi_0$ is an object $\xi\in \XX(A)$,
where $A$ is a complete local ring with residue field $K$,
and $\xi:\Spec A\to \XX$
which has the following lifting property: For any 2-commutative diagram
of solid arrows
\[
\xymatrix{
\Spec R/I \ar[d] \ar[r]& \Spec A \ar[d]^{\xi} \\
\Spec R \ar[r] \ar@{..>}[ur] & \XX \\
}
\]
where $R$ is an artin local ring with residue field $K$, and $I$ is a square
zero ideal of $R$,
there exists a dotted arrow filling in the diagram.

\begin{Proposition}
\label{existversalfinite}
Let $S$ be a locally noetherian scheme
and $\XX$ an algebraic stack locally of finite type over $S$.
Let $p$ be a closed point on $\XX$ in the sense of \cite[(5.2)]{LM}.
Then there exist a complete noetherian local $S$-ring $A$
and an object $\xi\in \XX(\Spec A)$ that has properties:
\begin{enumerate}
\renewcommand{\labelenumi}{(\theenumi)}

\item The residue field $K$ of $A$ is of finite type over $S$, and the restriction $\xi_K$ of $\xi$ to $K$
represents the point $p$.

\item $\xi$ is an effective versal deformation of $\xi_K$.

\end{enumerate}
\end{Proposition}

\Proof
Let $\YY_p$ denote the reduced closed substack
associated to the closed point $p$.
Let $f:X\to \XX$ be a smooth surjective morphism
from a scheme $X$.
Then $f^{-1}(\YY_p)$ is a closed subscheme of $X$.
Take
a closed point $x\in f^{-1}(\YY_p)$.
The residue field $K(x)$
is of finite type over $S$, and $\Spec K(x)\to \XX$
represents $p$.
Let $\hat{\OO}_{X,x}$ be the completion of the local ring at $x$.
Then the induced morphism $\Spec \hat{\OO}_{X,x}\to \XX$
is formally smooth at the closed point $x$.
Namely, $\Spec \hat{\OO}_{X,x}\to \XX$ has the lifting property
depicted above (see the definition of effective versal deformations).
Thus, it gives rise to a desired effective versal
deformation.
\QED

We shall refer to an object $\xi\in\XX(\Spec A)$ with properties (1), (2)
in Proposition~\ref{existversalfinite} as
an {\it effective versal deformation on} $A$ for $p$.

\begin{Lemma}
\label{algebraization}
Let $S$ be an excellent scheme
and $\XX$ an algebraic stack locally of finite type over $S$.
Let $p$ be a closed point on $\XX$.
Let $A$ be
a complete noetherian local $S$-ring with the maximal ideal $\mathfrak{m}$, and
let
an object $\xi\in \XX(\Spec A)$ be
an effective versal deformation for $p$.
Then there exist an affine $S$-scheme $U$, a closed point $u$ on $U$,
an isomorphism $\hat{\OO}_{U,u}\cong A$,
and a smooth morphism $\tilde{\xi}:U\to \XX$ such that
the restriction of $\tilde{\xi}$ to $\Spec \hat{\OO}_{U,u}$
is isomorphic to $\xi$ via $\hat{\OO}_{U,u}\cong A$.

\end{Lemma}

\Proof
In virtue of Artin's algebraization theorem (\cite{A1},  \cite[(10.10), (10.11)]{LM}, \cite{AD}),
there exist an affine $S$-scheme $U$, a closed point $u$ on $U$,
an isomorphism $\hat{\OO}_{U,u}\cong A$,
and a smooth morphism $\tilde{\xi}:U\to \XX$ such that
the system $\{\tilde{\xi}_n\}_{n\ge 0}$ is isomorphic to
$\{\xi_n\}_{n\ge 0}$ in $\lim\XX(A/\mathfrak{m}^{n+1})$.
Here $\tilde{\xi}_n$ (resp. $\xi_n$) denotes the restriction of $\tilde{\xi}$
(resp. $\xi$) to $\Spec A/\mathfrak{m}^{n+1}$.
Since $\XX$ is an algebraic stack, thus there exists
a natural equivalence $\XX(A)\cong \lim\XX(A/\mathfrak{m}^{n+1})$.
Therefore
the restriction of $\tilde{\xi}$ to $\Spec \hat{\OO}_{U,u}$
is isomorphic to $\xi$.
\QED

\begin{Definition}
\label{stablepoint}
Let $S$ be a excellent scheme
and $\XX$ an algebraic stack locally of finite type over $S$.
A closed point $p$ on $\XX$ is strongly p-stable if there exist a complete noetherian
local $S$-ring whose residue field is of finite type over $S$
and an effective versal deformation
$\xi\in \XX(\Spec A)$ for $p$, that has the following property ({\bf S}):
\begin{quotation}

There exists a flat normal closed subgroup
$\FF\subset \Aut_{\XX,A}(\xi)$ whose quotient $\Aut_{\XX,B} (\xi)/\FF$ is finite over $\Spec A$
and the following compatibility condition (C) holds:
Let $f,g:\Spec B \to \Spec A$ be two morphisms
of $S$-schemes such that their pullbacks $f^*\xi$, $g^*\xi$ are
isomorphic to each other.
Let $\Aut_{\XX,B} (f^*\xi)\to \Aut_{\XX,A} (\xi)$ and
$\Aut_{\XX,B} (g^*\xi)\to \Aut_{\XX,A} (\xi)$ be the induced morphisms,
and let $\FF_{f}\subset \Aut_{\XX,B} (f^*\xi)$ and $\FF_{g}\subset \Aut_{\XX,B} (g^*\xi)$
denote the pullbacks of $\FF$ respectively.
Then for any isomorphism $f^*\xi\to g^*\xi$
the induced isomorphism $\Aut_{\XX,B} (f^*\xi)\to\Aut_{\XX,B} (g^*\xi)$
gives rise to an isomorphism $\FF_f\to \FF_g$.
\end{quotation}

\end{Definition}

For the sake of simplicity, we shall refer to an effective versal deformation $\xi$ with property
({\bf S}) as a {\it strongly p-stable effective versal deformation}.

\begin{Remark}
\begin{enumerate}
\renewcommand{\labelenumi}{(\roman{enumi})}
\item To verify the condition (C) in Definition~\ref{stablepoint},
it suffices to check that there exists one isomorphism $f^*\xi\to g^*\xi$
inducing an isomorphism $\Aut_{\XX,B} (f^*\xi)\to\Aut_{\XX,B} (g^*\xi)$
that gives rise to $\FF_f\to \FF_g$.
Indeed, the difference of two isomorphisms $a,b:f^*\xi\rightrightarrows g^*\xi$
is $b^{-1}\circ a:f^*\xi\to f^*\xi$.
The automorphism $b^{-1}\circ a:f^*\xi\to f^*\xi$
gives rise to an inner automorphism $\Aut_{\XX,B} (f^*\xi)\to\Aut_{\XX,B} (f^*\xi)$
and it induces an isomorphism $\FF_f\to \FF_f$ because $\FF_f$ is
{\it normal} in $\Aut_{\XX,B} (f^*\xi)$.

\item We will use the condition (C) also when we define GIT-like stability.
Fortunately, in characteristic zero the condition (C) in GIT-like stability
is vacant.
\end{enumerate}
\end{Remark}

\begin{Remark}
Let $\mathcal{H}\to \Spec A$ be a {\it finite} scheme over a complete
local ring $A$.
Let $\hat{\mathcal{H}} \to \Spf A$ be the formal scheme
obtained by completion with respect to $\mathfrak{m}$-adic topology,
where $\mathfrak{m}$ is the maximal ideal of $A$.
Then $\hat{\mathcal{H}}$ is finite over $\Spf A$.
However, the converse is {\it not} true.
Namely, the finiteness of $\hat{\mathcal{H}} \to \Spf A$
does not imply that $\mathcal{H}\to \Spec A$ is finite.
To see this, assume that $A$ is a complete discrete valuation ring
with quotient field $K$.
Then $\Spec K \sqcup\Spec A\to \Spec A$ is not finite.
On the other hand, the associated formal scheme $\Spf A\to \Spf A$ is finite.
\end{Remark}

\begin{Remark}
As noted in Introduction, the letter ``p'' in the term ``p-stable'' stands for {\it pointwise} and
{\it potential}.
It seems natural that a pointwise stability is defined
by using the language of deformation theory.
In fact, the notion of effective versal deformations
plays ``the role of the completions of local rings'' on algebraic stacks.
Put another way, we have
Artin's criterion, which provides (only one) powerful and systematic
method for verifying algebraicity of stacks  (cf. \cite{A2}).
This criterion is described in terms of deformation theory.
Thus, our formulation fits in with Artin's representability
criterion. (However, note that in this paper any algebraic stack is assumed to
have quasi-compact and separated diagonal.)
Hence we make an effort to describe our stability in terms of deformation
theory.
\end{Remark}

\begin{Example}
To give a feeling for the stability introduced in Definition~\ref{stablepoint},
let us give some typical examples.
\begin{enumerate}
\renewcommand{\labelenumi}{(\theenumi)}

\item Every closed point on schemes and algebraic spaces
is strongly p-stable.
Let $\XX$ be an algebraic stack.
Suppose that for any affine scheme $\Spec A$
and any object $\xi\in \XX(\Spec A)$ the
automorphism group space $\Aut_{\XX,A} (\xi)$ is finite over $\Spec A$.
Then every closed point on $\XX$ is strongly p-stable.
These examples are in the realm of Keel-Mori theorem.

\item Let $X$ be an algebraic space
and $\XX$ an fppf gerbe over $X$, that is,
any point $x\in X$ admits an fppf neighborhood $U\to X$
such that the pullback $\XX\times_{X}U$ is
isomorphic to the classifying stack $B_U\mathcal{G}$ with some fppf group algebraic space $\mathcal{G}$ over $U$. 
Then every closed point
on $\XX$ is strongly p-stable.

\item Let $X$ be a {\it reduced} scheme of finite type and separated over the complex number field $\CC$.
Let $G$ be a reductive algebraic group over $\CC$, that acts on $X$.
Let $x$ be a closed point on $X$. Suppose that
$x$ is pre-stable in the sense of \cite[Definition 1.7]{GIT}.
Then the image of $x$ in the quotient stack $[X/G]$ is strongly
p-stable
(we will see this in section 6 in a more refined and generalized form).

\item Let $X$ be a projective scheme over a field $k$.
Let $G$ be an algebraic group over $k$. Let $\mathfrak{Bun}_{G}$ be a moduli stack
of principal $G$-bundles on $X$. Let $P$ be a $G$-bundle on $X$.
This stack $\mathfrak{Bun}_{G}$ is algebraic (see for example \cite[section 1]{AD}).
Note that for every $G$-bundle $\mathcal{P}$ over $T$, the automorphism group
scheme $\Aut_{\mathfrak{Bun}_{G},T}(\mathcal{P})$ contains $\textup{Cent}(G)\times_{k}T$, where $\textup{Cent}(G)$
 is the center of $G$.
Suppose that there exists an effective versal deformation $\xi\in \mathcal{M}(A)$ for $P$
such that $\Aut_{\mathfrak{Bun}_{G},A}(\xi)$ is finite modulo $\textup{Cent}(G)\times_{k}A$. Then the point corresponding to $P$ is strongly p-stable.

\end{enumerate}
\end{Example}

\section{Coarse moduli space (strongly p-stable case)}

The main purpose of this section is to prove

\begin{Theorem}
\label{main1}
Let $S$ be an excellent scheme.
Let $\XX$ be an algebraic stack locally of finite type over $S$.
Then strongly p-stable points on $\XX$ form an open substack $\XX^{st}$,
and there exists a coarse moduli map
\[
\pi:\XX^{st}\longrightarrow X,
\]
such that $X$ is an algebraic space locally of finite type over $S$.
Furthermore, the morphism $\pi$ is a quasi-finite and universally closed morphism and induces an isomorphism 
$\OO_{X}\to \pi_*\OO_{\XX}$.
If $X'\to X$ is a flat morphism of algebraic spaces, then
$\XX^{st}\times_{X}X'\to X'$ is also a coarse moduli map.
\end{Theorem}

\vspace{1mm}

Let $I\XX$ be the inertia stack of $\XX$.
For an $S$-scheme $Y$, $I\XX(Y)$ consists of pairs $(a, \sigma)$,
where $a\in \XX(Y)$, and $\sigma$ is an automorphism of $a$ in $\XX(Y)$.
A morphism $(a,\sigma)\to (a',\sigma')$ in $I\XX(Y)$
is a morphism $f:a\to a'$ in $\XX(Y)$ such that $\sigma'\circ f=f\circ\sigma$.
There is a natural forgetting representable morphism $I\XX\to \XX$ sending
$(a,\sigma)$ to $a$.
It is isomorphic to the fiber product $\XX\times_{\Delta,\XX\times_{S}\XX,\Delta}\XX$, where $\Delta:\XX\to\XX\times_{S}\XX$ is a diagonal morphism.
The forgetting morphism $I\XX\to \XX$
is isomorphic to the first projection
$\textup{pr}_1:\XX\times_{\Delta,\XX\times_{S}\XX,\Delta}\XX\to \XX$.
Thus, the morphism $I\XX \to \XX$ is of finite type and separated.

\vspace{2mm}

Let $P:U\to \XX$ be a smooth morphism from an affine $S$-scheme $U$. Let $\eta\in \XX(U)$ be the object corresponding to $P$.
Let $\Aut_{\XX,U} (\eta)\to U$ be the group algebraic space of
the automorphism group of $\eta$.
We have the natural isomorphism $\Aut_{\XX,U} (\eta)\to U\times_{\XX}I\XX$ over $U$.

Consider the following contravariant functor
\[
F:(\textup{affine}\ U\textup{-schemes})\to (\textup{Sets})
\]
which to any $f:Y\to U$ associates
the set of normal closed subgroup spaces $\FF\subset \Aut_{\XX,Y} (f^*\eta)$ over $Y$,
which have following properties (i), (ii), (iii):

\begin{enumerate}
\renewcommand{\labelenumi}{(\roman{enumi})}

\item $\FF$ is flat over $Y$.

\item $\Aut_{\XX,Y} (f^*\eta)/\FF$ is finite over $Y$.

\item If $\textup{pr}_1,\textup{pr}_2:\Aut_{\XX,Y}(f^*\eta)\times_{I\XX}\Aut_{\XX,Y}(f^*\eta)\rightrightarrows\Aut_{\XX,Y} (f^*\eta)$ are the natural projections,
then $\textup{pr}_1^{-1}(\FF)=\textup{pr}_2^{-1}(\FF)$.

\end{enumerate}

\begin{Proposition}
\label{localpresentation}
The functor $F$ is locally of finite presentation,
that is to say, for any inductive system
\[
A=\displaystyle\lim_{\longrightarrow} A_{\lambda}
\]
of finitely presented $U$-rings \textup{(}cf. \cite[Appendix A]{Milne}, \cite[section 1]{Art0}\textup{)}, the natural map
\[
\Pi:\displaystyle\lim_{\longleftarrow} F(\Spec A_{\lambda})\to F(\Spec \displaystyle\lim_{\longrightarrow} A_{\lambda})
\]
is an isomorphism.
\end{Proposition}

\Proof
It follows from the uniqueness part
of Theorem~\ref{limitclosed} that $\Pi$ is injective.
Next, we prove that $\Pi$ is surjective.
Let $f_{\lambda}:\Spec A_{\lambda}\to U$ be the morphism associated to a $U$-ring $A_{\lambda}$.
Let $\xi\in \XX(\Spec A)$ be the object corresponding to the composite $\Spec A\to U\to \XX$.
Let $\FF\to \Spec A$ be a closed group subspace in $\Aut_{\XX,A}(\xi)\to \Spec A$,
that has properties (i), (ii) and (iii).
By Theorem~\ref{limitclosed}, there exist
$\alpha\in I$ and a closed subspace $\FF_{\alpha}\subset \Aut_{\XX,A_{\alpha}} (f_{\alpha}^*\eta)$
which induces $\FF$ in $\Aut_{\XX,A} (\xi)$.
In addition, by Proposition~\ref{grouplimit} and Proposition~\ref{limit2}
we may assume that $\FF_{\alpha}$ is a normal subgroup space that is flat over $\Spec A_{\alpha}$.
Consider the quotient $\Aut_{\XX,A_{\alpha}}(f^*_{\alpha}\eta)/\FF_{\alpha}$.
It is an (quasi-separated) algebraic space of finite presentation over $S_{\alpha}$ by \cite[(10.13.1)]{LM}.
Since $\Aut_{\XX,A}(\xi)/\FF$ is finite over $\Spec A$,
by Proposition~\ref{limit2} there exists an arrow $\alpha\to \mu$
such that $\Aut_{\XX,A_{\mu}}(f^*_{\mu}\eta)/\FF_{\mu}$ is finite over $\Spec A_{\mu}$,
where $\FF_{\mu}=\FF_{\alpha}\times_{A_{\alpha}}A_{\mu}$.
This means that $\Pi$ is surjective.
\QED

\vspace{1mm}

{\it Proof of Theorem~\ref{main1}}.
Let $p$ be a strongly p-stable closed point.
Let $A$ be a complete
noetherian local ring $A$ whose residue field is of finite type over $S$,
and let $\xi\in \XX(\Spec A)$ be a strongly p-stable effective versal deformation on $A$ for $p$.
By Lemma~\ref{algebraization}, we may assume that there exist
an affine $S$-scheme $U$, a smooth morphism $P:U\to \XX$,
and a closed point $u$ such that $\hat{\OO}_{U,u}=A$
and the restriction of $P$ to $\Spec \hat{\OO}_{U,u}$ is $\xi$.

Let $F$ be the functor defined above.
Observe that $F(\Spec A)$ is nonempty.
Indeed, $\xi\in \XX(\Spec A)$ has the property ({\bf S}) in Definition~\ref{stablepoint}, thus there exists a flat normal closed subgroup $\FF\subset \Aut_{\XX,A} (\xi)$ with property ({\bf S}) in Definition~\ref{stablepoint}.
Clearly, it satisfies the above conditions (i), (ii).
It suffices to check (iii).
Note that for any affine $S$-scheme $V$,
the set $\Aut_{\XX,A}(\xi)\times_{I\XX}\Aut_{\XX,A}(\xi)(V)$
consists of quintuples
\[
(f,g:V\rightrightarrows \Spec A,\ \phi:f^*\xi\stackrel{\sim}{\to}f^*\xi,\ \psi:g^*\xi\stackrel{\sim}{\to}g^*\xi,\ \sigma:f^*\xi\stackrel{\sim}{\to}g^*\xi)
\]
where $f^*\xi$ and $g^*\xi$ are the
pullbacks of $\xi$ by $f$ and $g$ respectively and $\psi=\sigma\circ\phi\circ\sigma^{-1}$.
The first (resp. second) projection
$\textup{pr}_1,\textup{pr}_2:\Aut_{\XX,A}(\xi)\times_{I\XX}\Aut_{\XX,A}(\xi)(V)\to \Aut_{\XX,A}(\xi)(V)$
sends $(f,g,\phi,\psi,\sigma)$ to $(f,\phi)$ (resp. $(g,\psi)$).
Therefore $\textup{pr}_1^{-1}(\FF)$ (resp. $\textup{pr}_2^{-1}(\FF)$)
consists of quintuples
\[
(f,g:V\rightrightarrows \Spec A,\ \phi:f^*\xi\stackrel{\sim}{\to}f^*\xi,\ \psi:g^*\xi\stackrel{\sim}{\to}g^*\xi,\ \sigma:f^*\xi\stackrel{\sim}{\to}g^*\xi)
\]
such that $\phi\in \FF$ (resp. $\psi\in\FF$).
Thus to prove $\textup{pr}_1^{-1}(\FF)=\textup{pr}_2^{-1}(\FF)$
it suffices to prove the following: Let $f,g:V\rightrightarrows \Spec A$
be two morphisms,
and let $\phi:f^*\xi\stackrel{\sim}{\to}f^*\xi$ be
an automorphism which lies in $\FF(f:V\to \Spec A)$.
Assume that $f^*\xi$ is isomorphic to $g^*\xi$.
Then for any isomorphism $\sigma:f^*\xi\to g^*\xi$,
the composite $\sigma\circ \phi\circ \sigma^{-1}$ lies in $\FF(g:V\to \Spec A)$.
We can check it by using the condition (C) in Definition~\ref{stablepoint}.
Indeed, for an isomorphism $\sigma:f^*\xi\to g^*\xi$,
the induced isomorphism
\[
\Aut_{\XX,V}(f^*\xi)\to \Aut_{\XX,V}(g^*\xi)
\]
sends $\phi$ to $\sigma\circ \phi\circ \sigma^{-1}$.
Therefore the condition (C) in Definition~\ref{stablepoint} implies
$\textup{pr}_1^{-1}(\FF)=\textup{pr}_2^{-1}(\FF)$.
Hence $F(\Spec A)$ is nonempty.

By the algebraic approximation theorem \cite[(2.2)]{A}
and Proposition~\ref{localpresentation},
there exists an \'etale $U$-scheme $W\to U$
such that $F(W)$ is nonempty.
Then the composite $C:W\to U\to \XX$ is a smooth morphism.
Let $\iota\in \XX(W)$ be the object corresponding to $W\to \XX$.
Then there exists the natural cartesian diagram
\[
\xymatrix{
\Aut_{\XX,W}(\iota)\ar[r]^(0.65){h} \ar[d]& I\XX \ar[d]\\
W\ar[r] & \XX. \\
}
\]
Since $F(W)$ is nonempty, there exists a normal closed subgroup
$\FF_{\iota}\subset \Aut_{\XX,W}(\iota)$
with properties (i), (ii), (iii).
Let $\YY$ be the image of $W\to \XX$. It is an open substack of $\XX$.
Then the closed subspace $\FF_{\iota}$ descends to a closed substack $\FF\subset I\YY$
that is flat over $\YY$.
Here we claim that $\FF$ is a subgroup in $I\YY$,
more
precisely, $\FF$ is stable under the multiplication
$m:I\YY\times_{\YY}I\YY\to I\YY$, the inverse $i:I\YY\to I\YY$,
and the unit section $e:\YY\to I\YY$.
To see $m(\FF\times_{\YY}\FF)\subset \FF$,
consider the 2-commutative diagram
\[
\xymatrix{
\Aut_{\XX,W}(\iota)\times_{W}\Aut_{\XX,W}(\iota)\ar[r] \ar[d]^{(h,h)}& \Aut_{\XX,W}(\iota)\ar[d]^{h}\\
I\YY\times_{\YY}I\YY\ar[r]^{m} & I\YY \\
}
\]
where the upper horizontal arrow is the multiplication of $\Aut_{\XX,W}(\iota)$.
Since the inverse image $(h,h)^{-1}(\FF\times_{\YY}\FF)$ equals to
$\FF_{\iota}\times_{W}\FF_{\iota}$, the morphism
$\FF_{\iota}\times_{W}\FF_{\iota}\to \FF_{\iota}$ descends to
$\FF\times_{\YY}\FF\to \FF$.
In a similar way, we can easily see that $\FF$ is stable under
$i:I\YY\to I\YY$.
Since the unit section $W\to \Aut_{\XX,W}(\iota)$ factors through
$\FF_{\iota}$, thus the composite $L:W\to \Aut_{\XX,W}(\iota)\to I\YY$ factors through
$\FF$. Note that $L$ is isomorphic to $W\to\YY\stackrel{e}{\to} I\YY$.
Since $W\to \YY$ is essentially surjective, we conclude that
$e:\YY\to I\YY$ factors through $\FF$.
Hence $\FF$ is a subgroup of $I\YY$.

By Proposition~\ref{existversalfinite}, for any closed point $y$ on $\YY$ there exist
a complete noetherian local $S$-ring $B$ whose residue field is
of finite type over $S$
and an effective versal deformation $\zeta\in\YY(B)$ for $y$.
Then $\FF$ gives rise to a normal closed subgroup space 
in $\Aut_{\YY,B}(\zeta)$ with properties (i), (ii), (iii).
Thus every closed point on $\YY$
is strongly p-stable.

By Theorem~\ref{rigidification}, there
exists an fppf gerbe $\rho:\YY\to \YY'$
such that for any affine $S$-scheme $V$ and
any object $a\in\YY(V)$ the homomorphism
$\Aut_{\YY,V}(a)\to \Aut_{\YY',V}(\rho(a))$ is surjective and its kernel
is $\FF\times_{\YY,a}V\subset \Aut_{\YY,V}(a)$.
Since $\Aut_{\XX,V}(a)/(\FF\times_{\YY,a}V)$ is finite for any $a\in\XX(V)$,
thus the inertia stack $I\YY'$ is finite over $\YY'$.
By Keel-Mori theorem and Theorem~\ref{rigidification},
there exists a coarse moduli maps $\YY'\to Y'$, 
and the composite morphism $\YY\to\YY'\to Y'$ is also a coarse moduli map.
Also, note that $\YY\to \YY'$ is a quasi-finite and universally closed morphism
by (iii) in Theorem~\ref{rigidification}.
Thus, $\YY'\to Y$ is so since $\YY\to Y$ is proper and quasi-finite.
Hence for any strongly p-stable closed point $x$\ on $\XX$,
there exists an open substack $\UU\subset \XX$ containing $x$,
such that every closed point on $\UU$ is strongly p-stable, and
it has a coarse moduli map.
Therefore strongly p-stable points form an open substack $\XX^{s}$.
Using the universality of coarse moduli space,
we conclude that $\XX^{st}$
has a coarse moduli space $X$, that is
an algebraic space locally of finite type over $S$.
Moreover, the coarse moduli map $\pi:\XX^{st}\to X$ is quasi-finite and universally closed.
To see that $\OO_X\to \pi_*\OO_{\XX}$ is an isomorphism,
we may suppose that $X$ is an affine scheme and $\pi:\XX\to X$
is a composite morphism $\pi'\circ\rho:\XX\to\XX'\to X$, where $\XX\to \XX'$
is an fppf gerbe, $\XX'\to X$ is a coarse moduli map such that
$\OO_X\to \pi_*'\OO_{\XX'}$ is an isomorphism.
Then the natural morphism $\OO_{\XX'}\to\rho_*\OO_{\XX}$
induced by the fppf gerbe $\rho:\XX\to \XX'$ is an isomorphism.
Hence $\OO_X\to \pi_*\OO_{\XX}$ is an isomorphism.

Finally, we will prove that for any flat morphism $X'\to X$
of algebraic spaces, $\XX^{st}\times_{X}X'\to X'$ is also a coarse moduli space.
The claim is 
Zariski local on $X$. Thus, we may assume that $\XX^{st}\to X$ is divided
into
$\XX^{st}\to \XX^{st'}\to X$, where $\XX^{st}\to \XX^{st'}$ is an fppf gerbe
and the inertia stack of $\XX^{st'}$ is finite over $\XX^{s'}$.
Then $\XX^{st}\times_{X}X'\to \XX^{st'}\times_{X}X'$ is an fppf gerbe
and by Keel-Mori theorem $\XX^{st'}\times_{X}X'\to X'$ is a coarse moduli map.
Hence $\XX^{st}\times_{X}X'\to X'$ is also a coarse moduli space.
\QED

From the proof of Theorem~\ref{main1}, we also see:

\begin{Corollary}
\label{everystable}
If a closed point $p$ on $\XX$ is strongly p-stable,
then
every effective versal deformation for $p$ is a strongly p-stable effective versal deformation.
\end{Corollary}

We will denote by $\XX^{st}$
the open substack consisting of strongly p-stable points on $\XX$.

\begin{Remark}
In general, the coarse moduli map $\XX^{st}\to X$ in Theorem~\ref{main1}
is not separated. Let $\XX=BG$ be the classifying stack of an affine
group scheme $G$ over a field $k$.
Then $BG=BG^{st}$
and the structure morphism $BG\to\Spec k$ is a coarse moduli map.
Assume that $G$ is not proper. Then $BG$ is not separated over $k$.
\end{Remark}

\begin{Proposition}
\label{nonempty}
Suppose that $\XX$ has a reduced open substack.
Then $\XX^{st}$ is not empty.
\end{Proposition}

\Proof
By \cite[6.11]{SGA1}, there exists a nonempty open substack $\YY\subset \XX$
such that the natural projection
$I\XX\times_{\XX}\YY\to \YY$ is flat.
Then each closed point on $\YY$ is strongly p-stable.
\QED

\section{Deformation of coarse moduli spaces and p-stable points}

In this section, we give a construction of
deformations of coarse moduli spaces,
which is one of key ingredients to a construction of coarse moduli spaces.
Let
\[
\xymatrix{
\XX_0 \ar[r] \ar[d]^{\pi_0}&  \XX \\
X_0  \\
}
\]
be a diagram, where $\XX_0\to\XX$ is a nilpotent deformation of an algebraic stack $\XX_0$
and $\pi_0$ is a coarse moduli map which induces a natural isomorphism
$\OO_{X_0}\to \pi_*\OO_{\XX_0}$.
We want to construct
a coarse moduli space for $\XX$ by deforming $X_0$
(under a certain natural setting).
We refer to such a construction as the deformation of a coarse moduli space.
Before proceeding into detail,
it seems appropriate to begin by observing some motivating examples:

\begin{Example}
\label{gitliftexample}
Let $X=\Spec R$ be an affine scheme of finite type
over a field $k$ and $X_{\textup{red}}=\Spec R/I$ the associated reduced
scheme. Let $G$ be a linearly reductive
algebraic group over $k$.
Suppose that $G$ acts on $X$, and it gives rise to an action on $X_{\textup{red}}$.
Assume that these actions are closed. (An action is closed if
every orbit is closed after some extension of the base field
(cf. \cite[Definition 0.8]{GIT}).
Namely, every point on $X$ is stable in the sense of Geometric
Invariant Theory \cite{GIT}.
Then the geometric quotient for the action on $X$
is $\Spec R^G$. Here $R^G$ denotes the invariant ring.
On the other hand,
the geometric quotient for the action on $X_{\textup{red}}$
is $\Spec (R/I)^G$.
We can view $\Spec R^G$ and $\Spec (R/I)^G$
as ``coarse moduli schemes'' for the quotient stacks
$[X/G]$ and $[X_{\textup{red}}/G]$ respectively.
Since $G$ is linearly reductive,
the natural map $R^G\to (R/I)^G$ is surjective.
Therefore the nilpotent deformation
$[X_{\textup{red}}/G]\to [X/G]$ induces
a deformation $\Spec (R/I)^G\to \Spec R$
of the coarse moduli scheme $\Spec (R/I)^G$,
and the coarse moduli scheme for $[X/G]$ can be obtained by
deforming that of $[X_{\textup{red}}/G]$.
This is the reason why the stability in the Geometric Invariant Theory
makes reference only to the properties of underlying orbits.
\end{Example}

\begin{Example}
\label{unipotent}
Let $\mathbb{A}^2=\Spec k[x,y]$ be an affine space over a field $k$.
Let $G_{a}=\Spec k[t]$ be an additive group.
Consider the action of $G_a$ on $\mathcal{A}^2$,
described by $(x,y)\mapsto (x,y+t)$.
Set $X_0=\Spec k[x,y]/(x)\subset \mathbb{A}^2$ and $X_1=\Spec k[x,y]/(x^2)\subset \mathbb{A}^2$. Then the coarse moduli spaces of the quotient stacks
$[X_0/G_a]$ and $[X_1/G_a]$ are
\[
\Spec (k[x,y]/(x))^{G_a}=\Spec k
\]
and
\[
\Spec (k[x,y]/(x^2))^{G_a}=\Spec k[x]/(x^2)
\]
respectively.
In particular, $(k[x,y]/(x^2))^{G_a}\to (k[x,y]/(x))^{G_a}$ is surjective.
Next consider the action on $\mathbb{A}^2$ given by
$(x,y)\mapsto (x, tx+y)$.
In this case, the closed substack $[X_0/G_a]$ in $[\mathbb{A}^2/G_a]$
is isomorphic to $\Spec k[y]\times_k BG_a$
and its coarse moduli space is $\Spec k[y]$.
However, the natural map $(k[x,y]/(x^2))^{G_a}\to (k[x,y]/(x))^{G_a}\cong k[y]$
is {\it not surjective}.
\end{Example}

From the examples in Example~\ref{gitliftexample} and~\ref{unipotent},
to construct the deformation of a coarse moduli space,
we need to impose a lifting property on invariant rings.
The purpose of this section is to axiomatize the property in the
framework of algebraic stacks and construct the deformation of a coarse
moduli space.
Let $\XX$ be an algebraic stack over a scheme $S$.
Let $\XX_0$ be a closed substack of $\XX$, which is
determined by a nilpotent coherent ideal sheaf $\mathcal{I}$.
Suppose that there exists a coarse moduli map
$\pi_0:\XX_0\to X_0$ such that $\OO_{X_0}\to \pi_{0*}\OO_{\XX_0}$ is
isomorphism, and for any flat morphism $X_0'\to X_0$ of algebraic spaces
$\XX_0\times_{X_0}X_0'\to X_0$ is also a coarse moduli map.

\begin{Proposition}
\label{set1}
For any \'etale morphism $U_0\to X_0$ from a scheme $U_0$,
there exists an \'etale (representable) morphism
$\mathcal{U}\to \XX$ such that $\mathcal{U}\times_{\XX}\XX_0\cong U_0\times_{X_0}\XX_0$. That is, there exists an \'etale deformation of $\textup{pr}_2:U_0\times_{X_0}\XX_0\to \XX_0$
to $\XX$.
Moreover, such a deformation is unique up to unique isomorphism.
\end{Proposition}

\Proof
Without loss of generality, we may assume $\II^2=0$.
Let $\mathsf{L}_{(U_0\times_{X_0}\XX_0)/\XX_0}$ and $\mathsf{L}_{U_0/X_0}$
denote the cotangent complexes of $\textup{pr}_2:U_0\times_{X_0}\XX_0\to \XX_0$
and $U_0\to X_0$ respectively (cf. \cite[(17.3)]{LM}, \cite{OL}).
If $\textup{pr}_1:U_0\times_{X_0}\XX_0\to U_0$ is the first projection,
then by \cite[17.3 (4)]{LM}, we have
\[
\mathsf{L}_{(U_0\times_{X_0}\XX_0)/\XX_0}\cong L\textup{pr}_{1}^*\mathsf{L}_{U_0/X_0}\cong L\textup{pr}_{1}^*\Omega_{U_0/X_0}\cong 0.
\]
According to Olsson's deformation theory of representation morphisms
of algebraic stacks \cite[Theorem 1.4]{OL}, there exists
an obstruction $o\in \Ext^2(\mathsf{L}_{(U_0\times_{X_0}\XX_0)/\XX_0},\textup{pr}_2^*\II)$ whose vanishing is necessary and sufficient for the existence
of an \'etale deformation of $\textup{pr}_2$ to $\XX$.
Thus, there exists a desired \'etale morphism $\mathcal{U}\to \XX$
since
$\Ext^2(\mathsf{L}_{(U_0\times_{X_0}\XX_0)/\XX_0},\textup{pr}_2^*\II)=0$.
Furthermore, again by \cite[Theorem 1.4]{OL},
the vanishing $\Ext^0(\mathsf{L}_{(U_0\times_{X_0}\XX_0)/\XX_0},\textup{pr}_2^*\II)=\Ext^1(\mathsf{L}_{(U_0\times_{X_0}\XX_0)/\XX_0},\textup{pr}_2^*\II)=0$
imply the uniqueness of such a deformation.
\QED

In the above situation, we will refer to $\UU\to \XX$ (or simply $\UU$)
as an \'etale morphism associated to $U_0\to X_0$, and write
$\XX_U$ for $\UU$.

\begin{Lemma}
\label{set2}
Let $U_0\to X_0$ and $U_0'\to X_0$ be schemes that are \'etale over $X_0$. Let $U_0'\to U_0$ be a morphism over $X_0$.
Then there exists a morphism $\XX_{U'}\to \XX_{U}$ which makes
the diagram
\[
\xymatrix@R=4mm @C=17mm{
 \XX_{U_0'}\ar[rr] \ar[rd]\ar[dd]& & \XX_{U'} \ar[rd] \ar[dd]&  \\
  & \XX_0 \ar[rr] & & \XX \\
  \XX_{U_0} \ar[ru] \ar[rr]&  & \XX_U \ar[ru]& \\
}
\]
2-commutative, where $\XX_{U_0'}=U_0'\times_{X_0}\XX_0$ and $\XX_{U_0}=U_0\times_{X_0}\XX_0$.
Such a morphism is unique up to isomorphism.

\end{Lemma}

\Proof
Note first that $\XX_{U_0}\to U_0$
is a coarse moduli map.
Applying Proposition~\ref{set1} to $U_0'\to U_0$ we have
an \~etale morphism $\XX_{U'}'\to \XX_{U}$  associated to
the \'etale morphism $U_0'\to U_0$, that is,
a unique deformation of $\XX_{U_0'}\to \XX_{U_0}$.
Then the composite $\XX_{U_0'}'\to \XX_{U_0}\to \XX$
is also an \'etale morphism associated to $U_0'\to X_0$.
Finally, the uniqueness of $\XX_{U'}$ implies our assertion.
\QED


Consider the following lifting property $(\mathbb{L})$.

\vspace{1mm}
$(\mathbb{L})$: For any \'etale morphism $U_0\to X_0$ from an affine scheme $U_0$
the natural map $\Gamma(\XX_{U},\OO_{\XX_{U}})\to \Gamma(\XX_{U_0},\OO_{\XX_{U_0}})$ is surjective, where $\XX_{U_0}$ denotes $U_0\times_{X_0}\XX_0$.

\vspace{1mm}

If $(\mathbb{L})$ is satisfied, we say that $\XX$ has the property $(\mathbb{L})$ with respect to $\XX_0$.

\begin{Proposition}
\label{deformmoduli}
Assume that $\XX$ has the property $(\mathbb{L})$ with respect to $\XX_0$.
Then there exists a commutative diagram
\[
\xymatrix{
\XX_0 \ar[r] \ar[d]^{\pi_0} & \XX \ar[d]^{\pi} \\
X_0 \ar[r] & X, \\
}
\]
where $\pi$ is a coarse moduli map for $\XX$, which induces
a natural isomorphism
$\OO_X\to \pi_*\OO_{\XX}$.
Furthermore, $X_0\to X$ is a nilpotent deformation
and $X$ is locally of finite type.
For any flat morphism $X'\to X$ of algebraic spaces,
$\textup{pr}_2:\XX\times_XX'\to X'$ is also a coarse moduli map.
\end{Proposition}

The proof proceeds in several steps.

\begin{Lemma}
\label{step1}
Let $U_0'\to U_0$ be an \'etale morphism of affine schemes
over $X_0$.
Suppose that the natural map
$\Gamma(\XX_{U},\OO_{\XX_{U}})\to \Gamma(\XX_{U_0},\OO_{\XX_{U_0}})=\Gamma(U_0,\OO_{U_0})$ is surjective, where $\XX_{U_0}$ denotes $U_0\times_{X_0}\XX_0$. Let $U:=\Spec \Gamma(\XX_{U},\OO_{\XX_{U}})$.
Let $U'$ be a unique \'etale deformation of $U_0'\to U_0$ to $U$ $($cf. \cite[Theorem 3.23]{Milne}$)$. Let $\XX_{U'}'=\XX_U\times_{U}U'$.
Then there exists a unique (up to unique 2-isomorphism)
isomorphism $\XX_{U'}'\to \XX_{U'}$
of deformations of $\XX_{U_0'}:=\XX\times_{X_0}U_0'\to \XX_{U_0}$ to $\XX_U$.
Furthermore, there exists a natural isomorphism
$\Gamma(\XX_{U'},\OO_{\XX_{U'}})=\Gamma(U',\OO_{U'})$.
In particular, the natural map
$\Gamma(\XX_{U'},\OO_{\XX_{U'}})\to
\Gamma(\XX_{U_0'},\OO_{\XX_{U_0'}})=\Gamma(U_0',\OO_{U_0'})$ is surjective.
\end{Lemma}

\Proof
We prove the first claim.
Note that by Lemma~\ref{set2}
there exists $\XX_{U'}\to \XX_U$ that can be viewed as a unique \'etale deformation
of $\XX_{U_0'}\to \XX_{U_0}$ to $\XX_U$.
Thus by Proposition~\ref{set1} it suffices to prove that
$\XX_{U'}'\times_{\XX_U}\XX_{U_0}$ is isomorphic to $\XX_{U_0'}$
over $\XX_{U_0}$.
To see this, notice that there is a natural closed
immersion $\XX_{U_0}\to \XX_{U}\times_{U}U_0$ over $U_0$.
This morphism is not necessarily an isomorphism.
More generally, we have a diagram
\[
\xymatrix{
 \XX_{U_0'}\ar[r] \ar[d]& \XX_{U'}'\times_{U'}U_0' \ar[r] \ar[d] & U_0' \ar[d] \\
  \XX_{U_0} \ar[r] &\XX_{U}\times_{U}U_0 \ar[r] & U_0 \\
  }
\]
where all squares are cartesian, and $\XX_{U_0}\to \XX_{U}\times_{U}U_0$
and $\XX_{U_0'}\to \XX_{U'}'\times_{U'}U_0'$ are
closed immersions.
Since $\XX_{U_0}\times_{(\XX_{U}\times_{U}U_0)}(\XX_{U'}'\times_{U'}U_0')$
is naturally isomorphic to $\XX_{U_0}\times_{\XX_{U}}\XX_{U'}'$,
thus we have $\XX_{U_0}\times_{\XX_{U}}\XX_{U'}'\cong \XX_{U_0'}$
over $U_0'$. This implies the first claim.
Next we prove the second claim.
Let $p:Z\to \XX_U$ be a smooth surjective morphism from an affine scheme $Z$.
Then there exists an exact sequence
\[
\Gamma(U,\OO_{U})=\Gamma(\XX_{U},\OO_{\XX_{U}})\stackrel{p^*}{\to} \Gamma(Z,\OO_Z)\stackrel{\textup{pr}_{1,2}^*}{\rightrightarrows} \Gamma(Z\times_{\XX_{U}}Z,\OO_{Z\times_{\XX_U}Z}).
\]
Since $U'\to U$ is flat, thus the exactness holds after base changing
by $U'\to U$. This implies the second claim because $\XX_{U'}\cong\XX_{U'}'$.
\QED

\begin{Remark}
\label{flatinv}
By the above proof,
the property that $\OO_X\to \pi_*\OO_{\XX}$ is an isomorphism
is stable under flat base changes on $X$.
\end{Remark}

Let $X_{0,\textup{\'et}}$ be the \'etale site, whose objects
are affine schemes over $X_0$, and whose morphisms are $X_0$-morphisms.
Let $\OO$ be a sheaf of rings, which is determined by
\[
\OO(U_0):=\OO(U_0\to X_0)=\Gamma(\XX_{U},\OO_{\XX_U})
\]
for any \'etale morphism $U_0\to X_0$ in $X_{0,\textup{\'et}}$.
For any $U_0\to X_0$
the natural morphism $\OO(U_0\to X_0)\to\OO_{X_0}(U_0\to X_0)$
is surjective since $\XX$ has the property $(\mathbb{L})$ with respect
to $\XX_0$.
By Lemma~\ref{step1}, for any morphism $U_0'\to U_0$ in $X_{0,\textup{\'et}}$,
there exists a natural isomorphism
$\OO_{X_0}(U_0)\otimes_{\OO(U_0)}\OO(U_0')\cong \OO_{X_0}(U_0')$.
Thus, the sheaf $\OO$ induces a deformation of the ringed topos associated
to $(X_{0,\textup{\'et}},\OO_{X_0})$.
Namely, it gives rise to a deformation $X_0\hookrightarrow X$ of the algebraic space
$X_0$.

\begin{Lemma}
\label{step2}
There exists a morphism $\pi:\XX\to X$
that has the following properties:
\begin{enumerate}
\renewcommand{\labelenumi}{(\alph{enumi})}

\item The diagram
\[
\xymatrix{
 \XX_0 \ar[r]\ar[d]& \XX \ar[d] \\
  X_0 \ar[r] & X\\
  }
\]
is commutative;

\item the natural morphism $\OO_X\to \pi_*\OO_{\XX}$ is an isomorphism.
\end{enumerate}
\end{Lemma}

\begin{Remark}
The diagram in Lemma~\ref{step2} is not necessarily cartesian.
\end{Remark}

\Proof
Let $Z\to X$ be an \'etale surjective morphism from an affine scheme $Z$.
Let $\XX_{Z}\to \XX$ be an \'etale morphism associated to
$\textup{pr}_2:Z\times_{X}X_0\to X_0$.
Namely, there exists a natural morphism $\XX_{Z}\to Z$.
Let $W\to \XX_Z$ be a smooth surjective morphism from a scheme $W$.
To construct a desired morphism $\XX\to X$, it suffices to show that
there exists a morphism $[W\times_{\XX}W\rightrightarrows W]\to [Z\times_{X}Z\rightrightarrows Z]$ of groupoids (cf. \cite[(2.4.3)]{LM}).
To this end, it is enough to prove that
there exists a morphism
$\XX_Z\times_{\XX}\XX_Z\to Z\times_XZ$
which makes the diagram
\[
\xymatrix{
 \XX_{Z}\times_{\XX}\XX_Z\ar@<-0.5ex>[r] \ar@<0.5ex>[r] \ar[d] &  \XX_Z \ar[d] \\
 Z\times_XZ\ar@<-0.5ex>[r] \ar@<0.5ex>[r]& Z \\
  }
\]
commute.
If $Z_0\to X_0$ denotes the projection $Z\times_{X}X_0\to X_0$,
then $\textup{pr}_{i}:\XX_{Z}\times_{\XX}\XX_Z\to \XX_Z$
($i=1,2$) is an (\'etale) deformation of projection
$\textup{pr}_{i}:\XX_0\times_{X_0}(Z_0\times_{X_0}Z_0)\to \XX_{Z_0}=\XX_0\times_{X_0}Z_0$ ($i=1,2$) respectively.
Also, $\textup{pr}_{i}:Z\times_{X}Z\to Z$
($i=1,2$) is an (\'etale) deformation of projection
$\textup{pr}_{i}:Z_0\times_{X_0}Z_0\to Z_0$ ($i=1,2$) respectively.
Thus by Lemma~\ref{step1}, the projections
$\XX_Z\times_{Z,\textup{pr}_i}(Z\times_XZ)\to \XX_Z$ ($i=1,2$)
is naturally isomorphic to $\textup{pr}_{i}:\XX_{Z}\times_{\XX}\XX_Z\to \XX_Z$
($i=1,2$).
Hence we can obtain a desired morphism $\XX_Z\times_{\XX}\XX_Z\to Z\times_XZ$.
Finally, by the construction we have $\Gamma(\XX_{Z},\OO_{\XX_{Z}})=\Gamma(Z,\OO_Z)$. This means that $\OO_X\to \pi_*\OO_{\XX}$ is an isomorphism.
\QED

\begin{Lemma}
\label{step3}
The morphism $\pi$ in Lemma~\ref{step2} is a coarse moduli map for $\XX$.
For any flat morphism $X'\to X$ of algebraic spaces,
$\textup{pr}_2:\XX\times_XX'\to X'$ is also a coarse moduli map.
\end{Lemma}

\Proof
First of all, it is clear that for any algebraically closed field $K$
the morphism
$\pi$ induces a bijective map from the set of isomorphism
classes of $\XX(K)$ to $X(K)$ because $\pi_0$ does.
If $X'\to X$ is a flat morphism of algebraic spaces,
then $\textup{pr}_2:\XX_0\times_{X_0}(X'\times_{X}X_0)\to X'\times_{X}X_0$
is a coarse moduli map by our assumption.
On the other hand, the underlying map of $\textup{pr}_2:\XX\times_XX'\to X'$
can be identified with that of $\textup{pr}_2:\XX_0\times_{X_0}(X'\times_{X}X_0)\to X'\times_{X}X_0$.
Therefore, $\textup{pr}_2$ induces a bijective map from the set of isomorphism
classes of $\XX\times_XX'(K)$ to $X'(K)$ for any algebraically closed
field $K$.
Thus, it remains to prove that for any flat morphism
$X'\to X$, $\textup{pr}_2:\XX\times_XX'\to X'$ is universal among
morphisms to algebraic spaces.
Let $f:\XX\to W$ be a
morphism to an algebraic space $W$.
What we have to prove is that there exists a unique morphism
$\phi:X\to W$ such that $f=\phi\circ \pi$.
We first prove the case where $W$ is a scheme.

\begin{Claim}
\label{coarsescheme}
If $W$ is a scheme, 
there exists a unique morphism
$\phi:X\to W$ such that $f=\phi\circ \pi$.
For any flat morphism $X'\to X$ of algebraic spaces,
the projection $\textup{pr}_2:\XX\times_{X}X'\to X'$ is also universal
among morphisms to schemes.
\end{Claim}

{\it Proof of Claim.}
The composite $\XX_0\to \XX\to W$ induces
a morphism $X_0\to W$.
It is enough to show that if $U\to X$  is an \'etale morphism
from an affine scheme $U$, then there exists a unique morphism
$U\to W$ extending $U_0:=U\times_XX_0\to X_0\to W$, such that $\XX_U\cong \XX\times_XU\to \XX\to W$ is equal to $\XX_U\to U\to W$.
Since $U$ and $U_0$ have the same underlying topological space,
thus we have a map $U\to W$ as topological spaces.
Thus it suffices to construct a morphism of structure sheaves.
To this end, we may assume that $W$ is affine.
Since we have the morphism $\XX_U\to W$,
there exists a morphism $\Gamma(W,\OO_W)\to \Gamma(\XX_U,\OO_{\XX_U})=\Gamma(U,\OO_U)$.
Clearly, it gives rise to a desired morphism.
The uniqueness follows from the fact that $U\to W$ should arise from
$\Gamma(W,\OO_W)\to \Gamma(U,\OO_U)$ associated to $\XX_U\to \XX\to W$.
By Remark~\ref{flatinv}, for any flat morphism $X'\to X$,
$\textup{pr}_2:\XX\times_{X}X'\to X'$ induces a natural isomorphism
$\OO_{X'}\to \textup{pr}_{2*}\OO_{\XX\times_{X}X'}$.
Therefore, the above argument can be applied to
$\textup{pr}_2:\XX\times_{X}X'\to X'$.
\QED

Next we will show the case 
when $W$ is an algebraic space by reducing it to Claim~\ref{coarsescheme}.
This part is done by a more or less well-known argument.
But for the reader's convenience we will present a proof here.
Let $V\to \XX$ be a smooth surjective morphism from an affine
scheme $V$, which gives rise to
\[
R=V\times_{\XX} V\stackrel{\textup{pr}_{i}}{\rightrightarrows} V \to \XX.
\]
Let $X'\to X$ be a flat morphism of algebraic space.
We will prove that the sequence
\[
W(X')\to W(X'\times_XV) \rightrightarrows W(X'\times_XR)
\]
is exact. Let $\XX'=\XX\times_{X}X'$, $V'=X'\times_XV$ and $R'=X'\times_XR$.

We first prove that $W(X')\to W(V')$ is injective.
Let $Y \to W$ is an \'etale surjective morphism
from an affine scheme $Y$.
The fiber product $Y_1:=Y\times_WY$ is quasi-affine (cf. \cite[(1.3)]{LM}).
Let $\xi,\eta:X'\rightrightarrows W$ be two morphisms to $W$.
Assume that $\xi, \eta$ have the same image in $W(V')$.
It suffices to show $\xi=\eta$.
Take an \'etale surjective morphism $X''\to X'$ so that
there exist lifts $\xi_0:X''\to Y$ and $\eta_0:X''\to Y$
of $\xi$ and $\eta$ respectively.
Since $X''\to X'$ is \'etale and surjective, it is enough to show that
two composites $X''\stackrel{\xi_0,\eta_0}{\rightrightarrows}Y\to W$
coincide.
It is equivalent to proving that
$(\xi_0,\eta_0):X''\to Y\times Y$
factors through the image of $(\textup{pr}_1,\textup{pr}_2):Y_1\to Y\times Y$.
Put $R''=R\times_{X}X''$,
$\XX''=\XX'\times_{X'}X''$ and $V''=V'\times_{X'}X''$.
Consider the commutative diagram
\[
\xymatrix{
V''\ar[d]\ar[r] & \XX'' \ar[d] \ar[r] & X''\ar@<0.5ex>[r]^{\xi_0} \ar@<-0.5ex>[r]_{\eta_0} \ar[d] & Y \ar[d] \\
V' \ar[r] & \XX' \ar[r] & X' \ar@<0.5ex>[r]^{\xi} \ar@<-0.5ex>[r]_{\eta} & W, \\
  }
\]
where the left and middle squares are cartesian diagrams.
Notice that the morphism $V''\to Y\times Y$ induced by
$(\xi_0,\eta_0)$ factors through the image of $(\textup{pr}_1,\textup{pr}_2):Y_1\to Y\times Y$ since $Y\to W$ is \'etale surjective and
$\xi, \eta$ have the same image in $W(V')$.
Denote by $\alpha\in Y_1(V'')$ the image.
Since we have the morphism $\XX''\to Y\times Y$ induced by
$(\xi_0,\eta_0)$,
thus $\textup{pr}_1^*(\alpha)$ and $\textup{pr}_2^*(\alpha)$
coincide in $Y_1(R'')$.
Since $Y_1$ is a scheme, by Claim~\ref{coarsescheme} we see that
$\alpha$ comes from $Y_1(X'')$.
Therefore, we conclude that $W(X')\to W(V')$ is injective.

Next we will prove the middle exactness
in $W(X')\to W(V') \rightrightarrows W(R')$.
Let $h:V'\to W$ be a morphism.
Suppose that the two composites $h\circ \textup{pr}_1,h \circ \textup{pr}_2:R'\to W$ coincide. (Namely, it gives rise to a morphism
$\XX'\to W$.)
It suffices to show that $h$ arises from some $X'\to W$.
For ease of notation, we replace $\XX'$, $X'$, $V'$ and $R'$ by $\XX$, $X$, $V$, and $R$
respectively.
The morphism $\XX_0\to \XX\to W$ induces
$X_0\to W$.
Take an \'etale surjective morphism
$X_0'\to X_0$ from an affine scheme $X_0'$
so that there exists a lift $X_0'\to Y$ of $X_0\to W$.
By \cite[Theorem 3.23]{Milne},
there exists a unique \'etale deformation $X'\to X$ of $X_0'\to X_0$.
Then it is enough to construct $f:X'\to Y$
and $g:X'':=X'\times_{X}X'\to Y_1$ which makes
the diagrams ($\clubsuit$)
\[
\xymatrix{
 X''\ar@<0.5ex>[r]^{\textup{pr}_1} \ar@<-0.5ex>[r]_{\textup{pr}_2} \ar[d]^g & X' \ar[d]^f \\
Y_1 \ar@<0.5ex>[r]^{\textup{pr}_1} \ar@<-0.5ex>[r]_{\textup{pr}_2} & Y, \\
  }
\qquad
\xymatrix{
 V' \ar[r] \ar[d] & V \ar[d]^{h} \\
Y \ar[r] & W \\
  }
\]
commute, where $V':=V\times_{X}X'$.
We first construct $f:X'\to Y$.
To this aim, notice that $V_0':=X'_0\times_{X}V\cong X_0'\times_{X_0}V_0\to V_0:=X_0\times_XV$
factors through $\textup{pr}_2:Y\times_WV_0\to V_0$.
Again by \cite[Theorem 3.23]{Milne},
$V'\to V$ factors through $\textup{pr}_2:Y\times_WV\to V$.
Thus, we have $f':V'\to Y\times_WV\to Y$, which fits in with
the right diagram of $(\clubsuit)$.
Next we will observe that
$f'\circ \textup{pr}_1, f'\circ\textup{pr}_2:R'=R\times_XX'\rightrightarrows Y$
coincide.
Clearly, it holds after restricting to $R_0'=R\times_XX_0'$.
This implies that two morphisms $R_0'\rightrightarrows R_0\times_WY$ induced
by $R_0'\rightrightarrows Y$ coincide.
The morphisms $R'\rightrightarrows R\times_WY$ induced
by $R'\rightrightarrows Y$ are \'etale deformations of $R_0'\rightrightarrows R_0\times_WY$ respectively.
Thus by \cite[Theorem 3.23]{Milne},
two morphisms $R'\rightrightarrows R\times_WY$ coincide.
Hence $ f'\circ\textup{pr}_1$ and $ f'\circ\textup{pr}_2$
coincide, and it gives rise to $\XX':=\XX\times_{X}X'\to Y$.
This induces $\Gamma(Y,\OO_{Y})\to \Gamma(X',\OO_{X'})=\Gamma(\XX',\OO_{\XX'})$. Let $f:X'\to Y$ be the induced morphism.
Next we construct $g:X''\to Y_1$.
By the above construction, we have $\XX'\to Y$.
It gives rise to $a:\XX'\times_{\XX}\XX'\to Y\times_WY$ because the right diagram
of $(\clubsuit)$ commutes.
Note that by Lemma~\ref{step1} $(\XX'\times_{\XX}\XX')\times_{\XX}\XX_0$ is
naturally isomorphic to $\XX_0\times_{X_0}(X_0'\times_{X_0}X_0')$,
whose coarse moduli space is $X_0'\times_{X_0}X_0'$.
Also, notice $X_0'\times_{X_0}X_0'\cong (X'\times_XX')\times_{X}X_0$.
Applying the schematic case
to $\XX'\times_{\XX}\XX'\to X'\times_XX'$
there exists a unique morphism $g:X'\times_XX'\to Y\times_WY$
such that $a=g\circ \pi''$, where $\pi''$ is the natural morphism
$\XX'\times_{\XX}\XX'\to X'\times_XX'$.
By the construction of $f$ and $g$,
they make diagrams $(\clubsuit)$ commute.
Hence we obtain a desired morphism $\phi:X\to W$.
\QED

{\it Proof of Proposition~\ref{deformmoduli}}.
We now obtain our Proposition~\ref{deformmoduli}
from Lemma~\ref{step2} and Lemma~\ref{step3}. 
\QED

Let $\XX$ be an algebraic stack locally of finite type
over an excellent scheme $S$.
Let $p$ be a closed point on $\XX$.
We say that $p$ is p-stable if there exist a quasi-compact open substack $\UU\subset \XX$
containing $p$, and a closed substack $\UU_0\subset \UU$ determined by
a nilpotent ideal sheaf, that has the properties:
\begin{enumerate}
\renewcommand{\labelenumi}{(\roman{enumi})}
\item $p$ is strongly p-stable over $\UU_0$;

\item $\UU$ has the property $(\mathbb{L})$ with respect to $\UU_0$.
\end{enumerate}

According to Theorem~\ref{main1} and Proposition~\ref{deformmoduli},
we have:

\begin{Proposition}
Let $\XX$ be an algebraic stack locally of finite type
over an excellent scheme $S$.
Then p-stable points form an open substack $\XX^s$, and there exists
a coarse moduli space
\[
\pi:\XX^s\longrightarrow X,
\]
which induces a natural isomorphism $\OO_X\to \pi_*\OO_{\XX}$.
Moreover, it is universally closed and quasi-finite.
If $X'\to X$ is a flat morphism of algebraic spaces,
then $\XX^s\times_{X}X'\to X'$ is also
a coarse moduli space.
\end{Proposition}

\begin{Remark}
Ultimately, we shall only be interested in GIT-like p-stability
which will be introduced in the next section.
The reason, however, that we introduced the ad hoc notion of (strong)
p-stability is to prove some properties of GIT-like p-stable points.
In addition, we hope that the machinery of p-stability will
be useful in the other situations.

The property $(\mathbb{L})$ may seem to be hard to verify.
Nonetheless, in the next section we will show that
the property $(\mathbb{L})$ is satisfied in the GIT-like stable case
which we think of as a reasonable setting.

%
\end{Remark}

For the later use, we need:

\begin{Lemma}
\label{helpgit}
Let $\XX_0\to\XX$ be a closed immersion defined by a nilpotent ideal $\II$.
Suppose that there exists a coarse moduli map $\pi_0:\XX_0\to U_0$ to an affine scheme,
and for any flat morphism $U_0'\to U_0$, the projection
$q:\XX_{U_0'}:=\XX_0\times_{U_0}U_0'\to U_0'$ is also a coarse moduli map.
Let $U_0'\to U_0$ be an \'etale surjective
morphism of affine schemes.
Let $\XX_{U'}\to \XX$ be an \'etale morphism associated
to $U_0'\to U_0$, that is, a unique \'etale deformation $\XX_{U'}\to \XX$
of $\XX_{U_0'}\to \XX_0$.
Suppose that $\Gamma(\XX_{U'},\OO_{\XX_{U'}})\to \Gamma(\XX_{U_0'},\OO_{\XX_{U_0'}})$ is surjective.
Then $\Gamma(\XX,\OO_{\XX})\to \Gamma(\XX_0,\OO_{\XX_0})$
is surjective. (In particular, the property $(\mathbb{L})$ is \'etale local
property.)
\end{Lemma}

\Proof
By induction, we may and will assume that
$\II$ is a square zero ideal.
Let $\mathcal{J}:=\II\otimes_{\OO_{\XX_0}}\OO_{\XX_{U_0}}$.
Since $\Gamma(\XX_{U'},\OO_{\XX_{U'}})\to \Gamma(\XX_{U_0'},\OO_{\XX_{U_0'}})$ is surjective and $U_0'$ is affine, we have $\Gamma(U'_0,R^1q_*\mathcal{J})=0$.
On the other hand, $\Gamma(U_0,R^1\pi_*\II)\otimes_{\Gamma(U_0,\OO_{U_0})}\Gamma(U_0',\OO_{U_0'})=\Gamma(U'_0,R^1q_*\mathcal{J})$.
Since $U_0'\to U_0$ is \'etale and surjective,
we see $\Gamma(U_0,R^1\pi_*\II)=0$.
This means that
$\Gamma(\XX,\OO_{\XX})\to \Gamma(\XX_0,\OO_{\XX_0})$
is surjective.
\QED

\section{GIT-like p-stability}

In this section, we now introduce the notion of{\it GIT-like p-stable points}.
In this section, we will work over a perfect base field $k$.
Let $\XX$ be an algebraic stack locally of finite type over $k$.

\begin{Definition}
\label{defgit}
Let $p$ be a closed point on $\XX$.
The point $p$ is GIT-like p-stable if there exists an effective versal
deformation $\xi\in \XX(A)$,
which has the following properties:
\begin{enumerate}
\renewcommand{\labelenumi}{(\alph{enumi})}

\item  The special fiber of $\Aut_{\XX,A}(\xi)\to \Spec A$ is
linearly reductive, that is, $\Aut_{\XX,A}(\xi)\times_{\Spec A}\Spec A/\mathfrak{m}$
is a linearly reductive algebraic group over $\Spec A/\mathfrak{m}$,
where $\mathfrak{m}$ is the maximal ideal of $A$.

\item If $I$ denotes the ideal generated by nilpotent elements in $A$,
then there exists a normal subgroup scheme $\mathcal{F}$
in $\Aut_{\XX,A}(\xi)\times_A\Spec (A/I)$ which is
smooth and affine over $\Spec A/I$, and whose geometric fibers are connected.
Furthermore, the quotient $\Aut_{\XX,A}(\xi)\times_A\Spec (A/I)/\mathcal{F}$
is finite over $A/I$, and the compatibility condition $(C)$ as in Definition~\ref{stablepoint} holds
for $\mathcal{F}\subset\Aut_{\XX,A}(\xi)\times_A\Spec (A/I)$.

\end{enumerate}

\end{Definition}

\begin{Remark}
\label{gitdefremark}
\begin{enumerate}
\renewcommand{\labelenumi}{(\roman{enumi})}
\item Let $K$ be a field and $K'\supset K$ an extension of fields.
An algebraic group $G$ over $K$ is linearly reductive if and only if so is $G\times_{\Spec K}\Spec K'$ over $K'$.
Therefore, to verify (a) in Definition~\ref{defgit},
it is enough to show that there exist a field $K$ and a morphism
$v:\Spec K\to \XX$, such that $v$ represents the point $p$ and
$\Aut_{\XX,K}(v)$ is linearly reductive over $K$.

\item According to \cite[$VI_B$ Corollarie 4.4]{SGA3}, to
check that a group scheme $\mathcal{F}$ is smooth over a reduced scheme
$S$, it is enough to prove that every fiber is smooth
and $\mathcal{F}$ is Zariski locally equidimensional over $S$.

\item The characteristic zero case is simpler than the general case.
If $\mathcal{G}:=\Aut_{\XX,A}(\xi)\times_A\Spec (A/I)\to \Spec (A/I)$ is equidimensional,
then by \cite[$VI_B$ Corollaire 4.4]{SGA3} the identity component
$\mathcal{G}^{0}$ (cf. \cite[$VI_B$ D\'efinition 3.1]{SGA3})
is a smooth open normal subgroup, whose geometric fibers are connected.
By \cite[$VI_B$ Proposition 3.3]{SGA3}, 
the compatibility condition $(C)$
holds for $\mathcal{G}^{0}$.
Hence in characteristic zero, the property (b) in Definition~\ref{defgit} is
satisfied if and only if the following conditions hold:
\begin{enumerate}
\renewcommand{\labelenumi}{(\theenumi)}

\item $\mathcal{G}\to \Spec (A/I)$ is equidimensional,

\item $\mathcal{G}^{0}$ is affine
over $\Spec (A/I)$,

\item $\mathcal{G}/\mathcal{G}^{0}$
is finite over $\Spec (A/I)$.

\end{enumerate}

\item Consider the case when $\XX$ has finite inertia stack.
Suppose that all closed points on $\XX$ are GIT-like p-stable.
Then in characteristic zero, $\XX$ is a Deligne-Mumford stack (\cite{DM}),
and in positive characteristic
$\XX$ is a tame stack introduced by Abramovich, Olsson, and Vistoli (\cite{AOV}).
\end{enumerate}
\end{Remark}

\begin{Proposition}
\label{reducedgit}
Let $\XX_0$ be the reduced stack associated to $\XX$.
Let $p$ be a GIT-like p-stable point on $\XX$ (or equivalently $\XX_0$).
Then there exists an open substack $\YY_0\subset \XX_0$
containing $p$, which has a
 coarse moduli map
$\pi_0:\YY_0\to Y_0$ such that
it induces an isomorphism $\OO_{Y_0}\to \pi_{0*}\OO_{\YY_0}$,
and for any flat morphism $Y_0'\to Y_0$ the pullback
$\YY_0\times_{Y_0}Y_0'\to Y_0'$ is also a coarse moduli map.
\end{Proposition}

\Proof
Let $\xi:\Spec A\to \XX$ be an effective versal deformation
for $p$, which satisfies the properties (a) and (b) in Definition~\ref{defgit}.
Note that if $\mathcal{I}$ denotes the ideal generated by nilpotent
elements in $\OO_{\XX}$, then the ideal $I$
of nilpotent elements in $A$ is the pullback of $\mathcal{I}$
via $\xi:\Spec A\to \XX$ because $\XX$ is excellent.
Then $p$ is a strongly p-stable point on $\XX_0$.
Thus by Theorem~\ref{main1} we obtain our Proposition.
\QED

Let $\XX$ be an algebraic stack locally of finite type over the base field $k$.
According to Proposition~\ref{reducedgit},
we will denote by $\XX^{\gs}$ the open substack of GIT-like p-stable points.
The main purpose of the remainder of this section
is to prove Theorem~\ref{maingit}, that is, to show the existence of
a coarse moduli map for $\XX^{\gs}$ by applying the results developed in
section 4.

\begin{Lemma}
\label{liftgit}
Let $\XX$ be an algebraic stack locally of finite type over a perfect
field $k$.
Let $\XX_0$ be the reduced stack associated to $\XX$.
Let $\XX_0^{\gs}$ be the open substack of GIT-like p-stable points.
Let $\XX_0^{\gs}\to X_0$ be a coarse moduli map (cf. Proposition~\ref{reducedgit})
Let $U_0\to X_0$ be an \'etale morphism from an affine scheme $U_0$.
Let $\XX_U\to \XX^{\gs}$ be an \'etale morphism associated to $U_0\to X_0$.
(See Proposition~\ref{set1}.)
Suppose that $\XX_{U_0}=\XX_0^{\gs}\times_{X_0}U_0$ has the form
$[V_0/G]$, where $V_0$
is an affine $U_0$-scheme and $G$ is a linearly reductive
algebraic group over $k$, which acts on $V_0$ over $U_0$.
Then there exist a nilpotent deformation $V_0\to V$ and an
action of $G$ on $V$, which extends the action on $V_0$,
such that $\XX_U$ is isomorphic to $[V/G]$.
\end{Lemma}

\Proof
We may and will assume $\mathcal{I}^2=0$.
By our assumption,
there exists a morphism $f:\XX_{U_0}\to BG$ corresponding
to the principal $G$-bundle $V_0\to [V_0/G]\cong \XX_{U_0}$.
To prove our claim,
we will show that
there exists a dotted morphism
filling the diagram
\[
\xymatrix{
\XX_{U_0} \ar[r] \ar[d]^f&  \XX_U \ar@{..>}[ld]\\
BG.  &   \\
}
\]
Note that $\XX_{U_0}\to BG$ is representable because it arises from
the $G$-morphism $V_0\to \Spec k$.
According to \cite[Theorem 1.5]{OL}, an obstruction
for the existence of a dotted arrow lies in
$\Ext^1(Lf^*\mathsf{L}_{BG/\Spec k},\mathcal{I})=\Ext^1(f^*\mathsf{L}_{BG/\Spec k},\mathcal{I})$, where $\mathsf{L}_{BG/\Spec k}$ denotes the cotangent complex.
We claim
\[
\Ext^1(Lf^*\mathsf{L}_{BG/\Spec k},\mathcal{I})=0.
\]
To prove our claim, we first show that
$\mathsf{L}_{BG/\Spec k}$ is of perfect amplitude in $[0,1]$.
Consider the composition $\Spec k\stackrel{w}{\to} [\Spec k/G]=BG\to \Spec k$,
where $\pi$ is the natural projection.
Then we have a distinguished triangle
\[
Lw^*\LL_{BG/\Spec k}\to \LL_{\Spec k/\Spec k}\to \LL_{\Spec k/BG}\to Lw^*\LL_{BG/\Spec k}[1].
\]
Since $\LL_{\Spec k/\Spec k}=0$,
$\LL_{\Spec k/BG}\cong Lw^*\LL_{BG/\Spec k}[1]$ (the symbol $\cong$ means
the existence of a quasi-isomorphism).
The morphism $w$ is flat surjective of finite type, and thus
it suffices to show that $\LL_{\Spec k/BG}$ is of perfect amplitude in
$[-1,0]$.
To this end,
consider the flat base change $z:\Spec k\times_{BG}\Spec k\cong G\to \Spec k$
of $w:\Spec k\to BG$,
we have $z^*\LL_{\Spec k/BG}\cong \LL_{G/\Spec k}$.
Note that the morphism $G\to \Spec k$ is complete intersection in the
sense of \cite[(19.3.6)]{EGA1}, and
thus according to \cite[Ch. III (3.2.6)]{I}
the cotangent complex $\LL_{G/\Spec k}$ is of perfect
amplitude in $[-1,0]$. Therefore we deduce that $\LL_{\Spec k/BG}$ is
of perfect amplitude in $[-1,0]$.
(If $G$ is smooth, then $\LL_{G/\Spec k}$ is of perfect
amplitude in $[0]$.)
Hence we conclude that $\LL_{BG/\Spec k}$ is of perfect amplitude in $[0,1]$.
Since $f$ is flat, thus $Lf^*\LL_{BG/\Spec k}\cong f^*\LL_{BG/\Spec k}$ is
of perfect amplitude in $[0,1]$.
Thus we have $\mathcal{RH}om(f^*\LL_{BG/\Spec k},\mathcal{I})\cong
\mathcal{RH}om(f^*\LL_{BG/\Spec k},\OO_{\XX_{U_0}})\otimes^{\mathbb{L}}\mathcal{I}$, and $\mathcal{RH}om(f^*\LL_{BG/\Spec k},\OO_{\XX_{U_0}})$
is of perfect amplitude in $[-1,0]$.
Note that $[V_0/G]\to BG$ is an affine morphism.
In addition, the push forward with respect to $BG\to \Spec k$,
is an exact functor from the category of quasi-coherent sheaves on $BG$
to that of quasi-coherent sheaves on $\Spec k$.
Therefore the global section functor
on $[V_0/G]\cong \XX_{U_0}$ is exact.
Taking into account local-global spectral sequence for Ext groups,
we see that $\Ext^1(Lf^*\mathsf{L}_{BG/\Spec k},\mathcal{I})=0$
and there exists the desired arrow $\XX_U\to BG$.
The pullback of the natural projection $\Spec k\to BG$ by
$\XX_U\to BG$ gives rise to a principal $G$-bundle
$V:=\Spec k\times_{BG}\XX_U\to \XX_U$.
Notice that $V\to \XX_U$
is a flat deformation of $V_0\to [V_0/G]\cong \XX_{U_0}$
to $\XX_U$. Thus $V$ is affine.
This means that $\XX_U$ has the form $[V/G]$
where $V$ is affine, as desired.
\QED

\begin{Remark}
\label{liftremark}
\begin{enumerate}
\renewcommand{\labelenumi}{(\roman{enumi})}
\label{smoothcotangent}
\item From the above proof, we see that
if $G$ is a smooth algebraic group over $k$, then
$\mathsf{L}_{BG/\Spec k}$ is of perfect amplitude in $[1]$.

\item Also, the same argument shows the following: Let $A$ be a local
$k$-ring with residue field $k$ and maximal ideal $\mathfrak{m}$.
Let $G$ be a linearly reductive
algebraic group over $k$. Let $\hat{\mathcal{X}}\to \Spf A$
be a formal algebraic stack, that is, an inductive system of 
algebraic stacks $\mathcal{X}_n\to \Spec A_{n}$, where $A_n=A/\mathfrak{m}^{n+1}$. Suppose that $\mathcal{X}_0=BG=[\Spec k/G]$.
Then the system $\XX_0\hookrightarrow \XX_1\hookrightarrow \XX_2\hookrightarrow \cdots \XX_n\hookrightarrow$ has the form
\[
[\Spec k/G]\hookrightarrow [\Spec B_1/G]\hookrightarrow [\Spec B_2/G]\hookrightarrow \cdots [\Spec B_n/G]\hookrightarrow \cdots,
\]
where $B_i$ is an artin local $k$-ring over $A_i$ for any $i$.
(To see this, replace $f:\XX_{U_0}\to BG$ in the proof of Lemma~\ref{liftgit} by $\textup{Id}:\mathcal{X}_0=BG\to BG$ and apply the same argument.)
Furthermore, if $G$ is smooth, we have $\Ext^i(\mathsf{L}_{BG/k},(\mathfrak{m}^n/\mathfrak{m}^{n+1})|_{BG})=0$ for $i=0,1$, which deduces, by \cite[Theorem 1.5]{OL}, that the
system of $G$-rings $\{B_i\}_{i\ge0}$ is unique up to isomorphism.

\end{enumerate}
\end{Remark}

\begin{Proposition}
\label{reducedquot}
Let $\XX_0$ be the reduced stack associated to $\XX$.
Let $\XX_0^{\gs}$ be the open substack of GIT-like p-stable
points and $\pi_0:\XX_0^{\gs}\to X_0$ the coarse moduli map.
Then for any closed point $p$ on $X_0$ there exists an \'etale neighborhood
$U_0\to X_0$ and a closed point $u\in U_0$ lying over $p$, such that
(i) if $k(u)$ denotes the residue field of $u$, then $\Spec k(u)\to X_0$
extends to $\alpha:\Spec k(u)\to \XX_0$,
(ii) $U_0$ is an affine $k(u)$-scheme, and (iii) $\XX_{U_0}=\XX_0^{\gs}\times_{X_0}U_0$ has
the form $[V_0/G]$, where $V_0$ is finite over $U_0$,
and $G\to \Spec k(u)$ is the automorphism group of $\alpha$.
\end{Proposition}

The proof of this Proposition proceeds in several steps: Lemma~\ref{constgp}, \ref{sacfilt}, \ref{nbdconst}, Proposition~\ref{embed} and~\ref{reducedquotgit}.

\begin{Lemma}
\label{constgp}
Let $L$ be a field.
Let $\mathcal{H}\to \Spec A$ be a group scheme that is affine and smooth
over $A$, where $A$ is a complete noetherian local $L$-ring
with residue field $L$.
Suppose that the fiber of $\mathcal{H}$
over the closed point of $A$
is linearly reductive, and all geometric fibres of $\mathcal{H}\to Spec A$
are connected.
Then there exists an isomorphism $\mathcal{H}\to H\times_L\Spec A$
of group schemes over $\Spec A$.
\end{Lemma}

\Proof
Let $\hat{\mathcal{H}}\to \Spf A$ be the
formal group scheme associated to $\mathcal{H}\to \Spec A$,
which we can view as a smooth deformation of $H$ to $\Spf A$.
Note that $H$ is linearly reductive and thus higher group
cohomology groups are trivial. Thus by the deformation theory
of group schemes (cf. \cite[Expose III (3.7)]{SGA3}),
we see that there exists a unique deformation of $H$ to $\Spf A$,
that is, $H\hat{\times}_L\Spf A$, and thus there
exists an isomorphism between $\hat{\mathcal{H}}$ and
$H\hat{\times}_L\Spf A$.
Let $\mathcal{H}om_A(H\times_LA,\mathcal{H})$
(resp. $\mathcal{H}om_A(\mathcal{H},H\times_LA)$) be a functor
which to any $S\to \Spec A$ associates
the set of homomorphisms $H\times_LA\times_AS\to \mathcal{H}\times_AS$
(resp. $\mathcal{H}\times_AS \to H\times_LA\times_AS$)
of group schemes over $S$.
According to \cite[Expose XIX 2.6]{SGA3},
$\mathcal{H}$ is a reductive group over $\Spec A$ (cf. \cite[Expose XIX 2.7]{SGA3}).
Then by \cite[Expose XXIV 7.2.3]{SGA3},
$\mathcal{H}om_A(H\times_LA,\mathcal{H})$ and
$\mathcal{H}om_A(\mathcal{H},H\times_LA)$
are represented by schemes locally of finite type and separated over $A$.
By the above observation, there exist inductive systems
$\{\Spec A/\mathfrak{m}_A^{n+1}\to \mathcal{H}om_A(H\times_LA,\mathcal{H})\}_{n\ge 0}$, $\{\Spec A/\mathfrak{m}_A^{n+1}\to \mathcal{H}om_A(\mathcal{H},H\times_LA)\}_{n\ge 0}$ arising from the isomorphism
$\hat{\mathcal{H}}\cong H\hat{\times}_L\Spf A$.
Here $\mathfrak{m}_A$ is the maximal ideal of $A$.
Since $\mathcal{H}om_A(H\times_LA,\mathcal{H})$ and
$\mathcal{H}om_A(\mathcal{H},H\times_LA)$
are schemes,
these systems are uniquely extended to
$u:\Spec A\to \mathcal{H}om_A(H\times_LA,\mathcal{H})$
and $v:\Spec A\to \mathcal{H}om_A(\mathcal{H},H\times_LA)$.
Let us denote by $u:H\times_LA\to \mathcal{H}$ and
$v:\mathcal{H}\to H\times_LA$ the corresponding morphisms respectively
(here we abuse notion).
The uniqueness also implies that
$u\circ v$ is the identity morphism $\mathcal{H}$.
Similarly, $v\circ u$ is the identity morphism of $H\times_LA$.
This completes the proof.
\QED

\begin{Lemma}
\label{sacfilt}
Let $p$ be a GIT-like p-stable closed point.
There exists an open substack $\UU\subset \XX_0$
containing $p$, that has the property:
there exists a closed subgroup $\mathcal{F}\subset I\UU$
such that $\mathcal{F}$
is smooth and affine over $\UU$, geometric fibers of $\mathcal{F}\to \UU$
are connected,
and $I\UU/\mathcal{F}$ is finite over $\UU$.
\end{Lemma}

\Proof
This proof is parallel to Proposition~\ref{localpresentation}
and the first, second and third paragraph of the proof of
Theorem~\ref{main1}.
Let $\xi:\Spec A\to \XX$ be an effective versal deformation
for $p$, which satisfies the properties (a) and (b) in Definition~\ref{defgit}.
By Lemma~\ref{algebraization},
we extend the versal deformation $\xi$
to a smooth morphism $P:U\to \XX$ where $U$ is an affine scheme having a closed point $u$ such that $A\cong \hat{\OO}_{U,u}$ and $\xi\cong P|_{A}$.
Consider the following contravariant functor
$F:(\textup{affine}\ U\textup{-schemes})\to (\textup{Sets})$
which to any $f:Y\to U$ associates
the set of normal closed subgroup spaces $\mathcal{G}\subset \Aut_{\XX,Y} (f^*\eta)$ over $Y$
with following properties (i), (ii), (iii):

\begin{enumerate}
\renewcommand{\labelenumi}{(\roman{enumi})}

\item $\mathcal{G}$ is smooth and affine over $Y$, and geometric fibers are
connected,

\item $\Aut_{\XX,Y} (f^*\eta)/\mathcal{G}$ is finite over $Y$,

\item if $\textup{pr}_1,\textup{pr}_2:\Aut_{\XX,Y}(f^*\eta)\times_{I\XX}\Aut_{\XX,Y}(f^*\eta)\rightrightarrows\Aut_{\XX,Y} (f^*\eta)$ are the natural projections,
then $\textup{pr}_1^{-1}(\mathcal{G})=\textup{pr}_2^{-1}(\mathcal{G})$.
\end{enumerate}
As in the proof of Proposition~\ref{localpresentation},
using a standard limit argument (cf. Appendix)
we see that
the functor $F$ is locally of finite presentation.
By the approximation theorem, there exists an \'etale neighborhood
$V\to U$ of $u$ such that $F(V)$ is not empty.
As in the proof of Theorem~\ref{main1}, this implies that there exist an open substack $\UU\subset \XX$
containing $p$ and a closed subgroup $\mathcal{F}\subset I\UU$
that has the desired property.
\QED

(The first setup): For simplicity, we will replace $\XX_0^{\gs}$ and $X_0$ by $\XX$ and $X$
respectively.
For any $U\to X$, we write $\XX_U$ for $\XX\times_XU$.
Let $p:\Spec K\to \XX$ be a geometric point and 
$q:\Spec K\to X$ the composite of $p$ and the coarse moduli map,
where $K$ is a separable (algebraic) closure of $k$.
Shrinking $\XX$ by Lemma~\ref{sacfilt} we assume that
there is an closed subgroup $\mathcal{F}\subset I\XX$,
such that $I\XX/\mathcal{F}$ is finite,
$\mathcal{F}\to \XX$ is smooth and affine,
and all geometric fibers of $\mathcal{F}\to \XX$ are connected.
Fix an \'etale neighborhood $U\to X$ of $q:\Spec K\to X$,
where $U$ is an affine scheme.
There exists a rigidification
$\XX_U\to \YY_U$ associated to $\mathcal{F}$ (cf. Appendix~\ref{rigidification}).
Note that
$\YY_U$ has finite inertia stack and each stabilizer is linearly reductive.
Then according to \cite[Proposition 3.2]{AOV} (and its proof),
there exist an \'etale neighborhood $U'\to U$ and a closed point $u'\in U'$
(the image of $\Spec K\to U'$)
such that (i) $\Spec k(u')\to U'$ extends to $\alpha:\Spec k(u')\to \XX_{U'}$
and $U'$ is a $k(u')$-scheme,
and
(ii)
if $G$, $H$ and $\bar{G}$ denote $I\XX_{U'}\times_{\XX_{U'},\alpha}k(u')=\Aut_{k(u')}(\alpha)$,
$\mathcal{F}\times_{\XX_{U'},\alpha}k(u')$ and $G/H$ respectively,
then the base change
$\YY_{U'}=\YY_{U}\times_UU'$ has the form $[W/\bar{G}]$,
where $W$ is an affine scheme which is finite over $U'$
and has the trivial fiber $W\times_{U'}\Spec k(u')\cong \Spec k(u')$.
For ease of notation, we replace $U'$ and $u'$ by $U$ and $u$ respectively.
Let $L:=k(u)$ be the residue field of $u$.


\begin{Lemma}
\label{nbdconst}
There exists an \'etale neighborhood $V\to U$ of $q$,
such that (i) the morphism $W_V:=W\times_{U}V\to \YY_U\times_UV$ lifts to
$W_V\to \XX_V$, and (ii)
the group scheme
$W_{V}\times_{\XX_{V}}\mathcal{F}\to W_{V}$
is isomorphic to the constant group scheme
$H\times_{L}W_{V}\to W_{V}$.

\end{Lemma}

\Proof
Clearly, we may assume that $\XX$ is quasi-compact.
Observe first that
it is enough to prove the case when the base field is algebraically closed
(in particular, $L=K$).
Note that any algebraic extension of $k$ is separable.
In addition, $X$ and $\XX$ is of finite presentation over $k$.
Thus, standard limit arguments show
that an \'etale morphism $V'\to U\times_LK$
with the desired properties always arises from some \'etale morphism
$V\to U$ with such properties. Therefore we may and will assume that $k=L=K$.

Next we prove that there exists an \'etale neighborhood $V\to U$
that satisfies (i).
Let $O$ be a strict henselization of the local ring
of $\OO_{U,u}$.
Since strict henselization commutes with finite extensions of rings,
thus $W':=W\times_{U}\Spec O$ is the disjoint union of
spectrums of strict
henselian local rings, that is, $W'\cong \amalg_{i=1}^{r}\Spec A_i$,
where $A_i$ is a strict
henselian local ring with residue field $K$ for all $i$.
Since $u\in U$ has the trivial fiber $\Spec K$, we
see $i=1$. We let $W':=\Spec A$.
The rigidifying morphism $\XX_U\to \YY_U$ is smooth, and thus there
exists a lifting $W'\to \XX_U$ of $W'\to \YY_U$.
We can write $\Spec O\to U$ as
$\Spec (\textup{colim}_{\lambda}B_{\lambda})$,
where $U_{\lambda}:=\Spec B_{\lambda}\to U$ are \'etale neighborhoods of the (geometric) point $u:\Spec K\to U$.
Note that $\XX_U$ is of finite type over $\YY_U$.
Hence there exists some \'etale neighborhood
$U_{\mu}\to U$
such that $W_{\mu}:=W\times_UU_{\mu}\to \YY_U$ extends to
$W_{\mu}\to \XX_U$.

Next we will show that there exists an \'etale neighborhood
$U_{\mu'}\to U_{\mu}$
that satisfies (ii).
Let $W^{\diamond}$ denote the spectrum of the completion $\hat{A}$
of $A$ with respect to the maximal ideal. In other words,
we also have $W^{\diamond}=\Spec A\hat{\otimes}_{O}\hat{O}$, where $\hat{O}$ is the completion
of $O$ with respect to the maximal ideal of $O$.
Then by Lemma~\ref{constgp}, the group scheme $W^{\diamond}\times_{\XX}\mathcal{F}$
is isomorphic to $W^{\diamond}\times_{K}H$ over $W^{\diamond}$.
Then considering the category of rings over $A$
and applying Artin's approximation
theorem,
we see that there exists
an isomorphism of group schemes between
$W'\times_{\XX_{U_{\mu}}}\mathcal{F}$ and
$H\times_KW'$
over $W'$.
Applying
Theorem~\ref{limit4} and~\ref{homlimit}
to the system $\{U_{\lambda'}\}_{\lambda'\to \lambda}$,
we conclude, by standard limit arguments, that there exists
$U_{\mu'}\to U_{\mu}$ such that
the group scheme $W_{\mu'}\times_{\XX_{U_{\mu'}}}\mathcal{F}\to W_{\mu'}$ is isomorphic to
$H\times_KW_{\mu'}\to W_{\mu'}$.
\QED

Since $\YY_U\cong [W/\bar{G}]$, we see that
$W_{U}\times_{\XX}I\XX/(W_{U}\times_{\XX}\mathcal{F})$
is embedded in $W_{U}\times_L\bar{G}$.
Let $G'$ be the reduced scheme associated to $G$
and set $\bar{G}':=G'/H$. The schemes $G'$ and
$\bar{G}'$ have naturally (smooth) group structures because
$L$ is perfect.
Next we prove:

\begin{Proposition}
\label{embed}
There exists an \'etale neighborhood $V\to U$ of $q$
such that the group scheme $(W_{V}\times_{\XX}I\XX)\times_{(W_{V}\times_L\bar{G})}(W_{V}\times_L\bar{G}')\to W_{V}$
can be embedded into $W_{V}\times_LG'$ as the inverse image of 
$(W_{V}\times_{\YY_U}I\YY_U)\times_{(W_{V}\times_L\bar{G})}(W_{V}\times_L\bar{G}')\subset W_{V}\times_L\bar{G}'$ 
under $W_{V}\times_L{G}'\to W_{V}\times_L\bar{G}'$.
\end{Proposition}

\Proof
\noindent
(Step 0)
As in the proof of Lemma~\ref{nbdconst}, we may assume that
the base field is algebraically closed. Thus we will let $k=L=K$.
Moreover,
by Lemma~\ref{nbdconst} we can take an \'etale neighborhood $V\to U$ that 
has the properties
(i) and (ii) in Lemma~\ref{nbdconst}.
For ease of notation, we may replace $V$ by $U$.

\vspace{2mm}

\noindent
(Step 1)
First we observe that it suffices to construct 
the desired embedding over $W^{\diamond}$.
(For the notation $W^{\diamond}$, see the proof of Lemma~\ref{nbdconst}.
We will continue to use notation in the proof of Lemma~\ref{nbdconst}.)
Assume that
the group scheme $\mathsf{G}:=(W^{\diamond}\times_{\XX}I\XX)\times_{(W^{\diamond}\times_K\bar{G})}(W^{\diamond}\times_K\bar{G}')\to  W^{\diamond}$
can be embedded into $W^{\diamond}\times_KG'$.
Consider the natural projection $W'\times_K{G'}\to W'\times_K\bar{G}'$
and take the inverse image $P$ of $W'\times_{\XX}I\XX/(W'\times_{\XX}\mathcal{F})\subset W'\times_K\bar{G}$ in $W'\times_K{G}'$.
Let $F$ be the functor which to any $a:Z\to \Spec O$
associates the set of isomorphisms of the group
schemes from $({W}_Z\times_{\XX}I\XX)\times_{(W_Z\times_K\bar{G})}(W_Z\times_K\bar{G}')$ to $a^*P$.
Using Theorem~\ref{limit4} and~\ref{homlimit}
we easily see that $F$ is locally of finite presentation.
Then we can apply Artin's approximation to conclude that
the group scheme $({W}'\times_{\XX}I\XX)\times_{({W}'\times_K\bar{G})}({W}'\times_K\bar{G}')$
can be embedded into $W'\times_KG'$.
By standard limit arguments, we see that
there exists an \'etale neighborhood $V\to U$ of $q$
such that the group scheme $(W_{V}\times_{\XX}I\XX)\times_{(W_{V}\times_K\bar{G})}(W_{V}\times_K\bar{G}')\to W_{V}$
can be embedded into $W_{V}\times_KG'$ in the desired way.

\vspace{2mm}

\noindent
(Step 2)
Next we prove that there exists an embedding of 
the scheme $\mathsf{G}$ into 
$W^{\diamond}\times_KG'$, that is, an isomorphism
between $\mathsf{G}$ and $P^{\diamond}:=P\times_{W'}W^{\diamond}$ as schemes.
Note that $W^{\diamond}\times_{\XX}\mathcal{F}$
is isomorphic to $W^{\diamond}\times_KH$ as group schemes.
The quotient $\mathsf{Q}:=\mathsf{G}/(W^{\diamond}\times_{\XX}\mathcal{F})$
is a (non-flat) finite group scheme over $W^{\diamond}$.
We have $\YY_U\cong [W_U/\bar{G}]$ and thus
$\mathsf{Q}$
is a closed subgroup scheme of $W^{\diamond}\times_K\bar{G}'$.
Let $\{\textup{Id}=g_1,g_2,\ldots,g_n\}$ be the set of ($K$-valued) points 
on $\bar{G}'$. (Note that $\bar{G}'$ is finite \'etale over $K$.
Namely, $\bar{G}'$ can be viewed as the finite group $\{\textup{Id}=g_1,g_2,\ldots,g_n\}$.)
The connected component of $\mathsf{Q}$ on which $g_i$ lies, denoted by $\mathsf{Q}_i$,
is isomorphic to the pullback
of the diagonal $W^{\diamond}\to W^{\diamond}\times_KW^{\diamond}$
by $(\textup{Id}_{W^{\diamond}},g_i):W^{\diamond}\to W^{\diamond}\times_KW^{\diamond}$.
Therefore, each connected component of $\mathsf{Q}$ can be
identified with a closed subscheme in $W^{\diamond}$.
Let us identify
$\mathsf{Q}_i$
with $M_i$, where $M_i$ is a closed subscheme
of $W^{\diamond}$.
Let $b:\mathsf{G}\to \mathsf{Q}$ be the projection.
Note that this projection is a principal $H$-bundle, and $M_i$
is the disjoint union of the spectrums of complete local rings.
($H$ is the identity component of $G'$.)
Also, the principal bundle of a smooth algebraic group over a strict henselian local ring is trivial.
Therefore $\mathsf{G}$ is isomorphic to
$H\times_K\mathsf{Q}$ 
as a scheme over $W^{\diamond}$.
We will regard $\mathsf{G}$ as $H\times_K\mathsf{Q}$.
(Here we do not take care of group structures.)
Thus, it is enough to prove that the group scheme
$\mathsf{G}\to W^{\diamond}$
is isomorphic to the group scheme
$H\times_K\mathsf{Q}\to W^{\diamond}$
(here $H\times_K\mathsf{Q}$ equips with the group structure arising from that of $G'$
in the natural way).

\vspace{1mm}

To this end, we fix some notation.
We write $\hat{W}'$ for the formal scheme associated to $W^{\diamond}$
and the adic topology on $\hat{A}$ arising from the maximal ideal.
Let $\hat{\mathsf{G}}$ (resp. $(H\times_K\mathsf{Q})^{\wedge}$)
be the formal scheme obtained from $\mathsf{G}\to W^{\diamond}$ (resp. $H\times_K\mathsf{Q}\to W^{\diamond}$)
by completion over $\hat{W}'$.
Let
$\hat{\mathsf{Q}}_{i}$ be the formal scheme obtained by completing
along the closed point of $\mathsf{Q}_i$.
Applying Remark~\ref{liftremark} (ii) to the base change $\hat{\XX}\to \Spf \hat{O}$ of $\XX\to X$ to $\Spf \hat{O}$, we see that $\hat{\mathsf{G}}$
is a closed formal group subscheme of $G\hat{\times}_K\hat{W}'$.
(The special fiber of $\hat{\XX}$ is a classifying stack $BG$
(cf. \cite[(11.3)]{LM}).)

\vspace{2mm}

\noindent
(Step 3)
Next we will show that there exists a closed immersion of formal group schemes
\[
\phi:\hat{\mathsf{G}}\to  G'\hat{\times}_K\hat{W}'
\]
over $\hat{W}^{'}$,
which identifies $\hat{\mathsf{G}}$ with $(H\times_K\mathsf{Q})^{\wedge}$.
From (Step 2), it is clear that there exists a closed immersion
of formal group schemes $\phi':\hat{\mathsf{G}}\to G\hat{\times}_K\hat{W}'$.
Moreover, we know the existence of an isomorphism
$\hat{\mathsf{G}}\cong (H\times_K\mathsf{Q})^{\wedge}$
of formal schemes.
Thus, it is enough to show that this immersion factors through
$G'\hat{\times}_K\hat{W}'\subset G\hat{\times}_K\hat{W}'$.
To this end,
observe first that the ``identity component'' of $\hat{\mathsf{G}}$
is the closed subscheme $H\hat{\times}_K\hat{W}'$ (in $G\hat{\times}_K\hat{W}'$).
The fiber of $\hat{\mathsf{G}}\to \hat{W}'$ over the closed point
is $G'$ and its identity component is $H$. Hence
it is enough to prove that any smooth deformation of $H$
to $\hat{W}^{'}$ that is embedded in $G^0\hat{\times}_K\hat{W}^{'}$,
is $H\hat{\times}_K\hat{W}^{'}$. Here $G^0$ is the identity component
of $G$.
In characteristic zero, we have $G^0=H$, thus our assertion is clear
since ``identity component'' of $(H\times_K\mathsf{Q})^{\wedge}$
is a constant deformation of $H$ to $\hat{W}'$, and any surjective endmorphism of a noetherian local ring is an
isomorphism.
In positive characteristic, note that by Nagata's classification of
linearly reductive groups (see for example \cite[page 27]{GIT}), $H$ is a torus.
Moreover, $G
^0/H$ is a linearly reductive group, thus by the classification in
\cite[Proposition 2.13]{AOV} we see that $G^0/H$ is a diagonalizable group.
Let $\hat{\mathcal{H}}$ be a smooth deformation of $H$
to $\hat{W}^{'}$ that is embedded in $G^0\hat{\times}_K\hat{W}^{'}$.
($\hat{\mathcal{H}}$ is a constant deformation of $H$.)
Consider the composite homomorphism $\hat{\mathcal{H}}\to G^0\hat{\times}_K\hat{W}^{'}\to (G^0/H)\hat{\times}_K\hat{W}^{'}$, where the second homomorphism
is the natural projection.
Since $\hat{\mathcal{H}}$ is isomorphic to $H\hat{\times}_K\hat{W}^{'}$
as formal group schemes, the composite 
$\hat{\mathcal{H}}\to (G^0/H)\hat{\times}_K\hat{W}^{'}$
comes
from a homomorphism of abelian groups.
Namely, if $H=\Spec K[M]$ and $G^0/H=\Spec K[N]$, then
the composite arises from a homomorphism $N\to M$.
However, any $N\to M$ is the trivial homomorphism since $M$ is free and $N\otimes_{\ZZ}\QQ=0$.
This means $\hat{\mathcal{H}}\subset H\hat{\times}_K\hat{W}^{'}$.
As in characteristic zero case, we deduce that
$\hat{\mathcal{H}}$ is equal to $H\hat{\times}_K\hat{W}^{'}$.
Using this, we will show that $\phi'$ factors through $G'\hat{\times}_K\hat{W}'\subset G\hat{\times}_K\hat{W}'$.
To see this, we may consider the problem
by restricting $\hat{\mathsf{G}}$
to some open and closed subgroup of
$\hat{\mathsf{G}}\hat{\times}_{\hat{W}'}\hat{\mathsf{Q}}_i$,
that is smooth over $\hat{\mathsf{Q}}_i$ for each $i$.
($\hat{\mathsf{G}}$ is the union of such subformal schemes.)
Namely, we may assume that $\hat{\mathsf{G}}$ is smooth
over $\hat{W}'$, that is, we have an isomorphism $\hat{\mathsf{G}}\cong
G'\hat{\times}_K\hat{W}'$ of formal group schemes.
Consider the composite $\rho\circ\phi':\hat{\mathsf{G}}\to G\hat{\times}_K\hat{W}'\to (G/G')\hat{\times}_K\hat{W}'$, where $\rho$ is the natural projection.
It suffices to prove that the image of $\rho\circ\phi'$ is trivial.
In characteristic zero case, it is clear.
In positive characteristic, we put $G/G'=\Spec K[N]$, where $N$ is an
abelian group whose order is a power of $\textup{ch}(k)$.
For each connected component $\hat{\mathsf{C}}_i$ of $\hat{\mathsf{G}}$, choose
a section $s_i:\hat{W}'\to \hat{\mathsf{C}}_i$.
Then the composite $\rho\circ\phi'$ is uniquely determined by the images $\rho\circ\phi'(s_i)$. Note that the number of connected components of $\hat{\mathsf{G}}$ is prime to $\textup{ch}(k)$.
Moreover, the identity component of $\hat{\mathsf{G}}$ maps to
the unit of $(G/G')\hat{\times}_K\hat{W}'$, hence for each $i$,
there exists a positive integer $n$, such that $\rho\circ\phi'(s_i^n)$ is the unit of $(G/G')\hat{\times}_K\hat{W}'$ and $n$ is prime to $\textup{ch}(k)$.
Therefore, we see that $\rho\circ\phi'(s_i)$
is the unit. This implies our claim.
Hence we have the desired morphism $\phi$.

\vspace{2mm}

\noindent
(Step 4)
We then claim that $\phi$ is induced by an isomorphism
$\mathsf{G}\to H\times_K\mathsf{Q} $ of group schemes
over $W^{\diamond}$.
For any $i$, let $\bar{G}'_i$ be a subgroup of $\bar{G}'$, which is generated by $g_i$. (Here we regard $\bar{G}'$ as a finite group.)
Let $G_i'$ be the preimage of $\bar{G}_i'$ under $G'\to \bar{G}'$.
For any $i$, set $\mathsf{G}_i:=G_i'\times_K\mathsf{Q}_i$  via a fixed
isomorphism $\mathsf{G}\cong H\times_K\mathsf{Q}$ of schemes. The scheme
$\mathsf{G}_i$ can naturally be viewed as a closed and open subgroup
scheme of $\mathsf{G}\times_{W^{\diamond}}\mathsf{Q}_i$ (over $\mathsf{Q}_i$).
(Note that at this point we do not claim that
the group structure of $\mathsf{G}_i$ comes from that of $G'_i$.)
On the other hand, we denote by $G'_i\hat{\times}_{K}\hat{\mathsf{Q}}_i$
the formal group scheme over $\hat{\mathsf{Q}}_i\subset W'$, whose group structure arises from
that of $G_i'$.
The formal group scheme $G'_i\hat{\times}_{K}\hat{\mathsf{Q}}_i\to \hat{\mathsf{Q}}_i$
is a closed and open subgroup scheme of $(H\times_K\mathsf{Q})^{\wedge}\hat{\times}_{\hat{W}'}\hat{\mathsf{Q}}_i$.
The morphism $\phi$ identifies the formal completion $\hat{\mathsf{G}}_i$
of $\mathsf{G}_i$
over $\hat{\mathsf{Q}}_i$ with
$G'_i\hat{\times}_{K}\hat{\mathsf{Q}}_i$ (as a formal group scheme).
We denote by $\phi_i$ the restriction of $\phi$ to $\hat{\mathsf{G}}_i$.
We will extend $\phi_i$ to an isomorphism $\mathsf{G}_i\to G'_i\times_K\mathsf{Q}_i$ of group schemes over $\mathsf{Q}_i$. (Note that the target $G'_i\times_K\mathsf{Q}_i$
has the group structure coming from $G_i'$.)
To this aim, consider the functor $\Hom^{\flat}_{\mathsf{Q}_i}(\mathsf{G}_i,G_i'\times_K\mathsf{Q}_i)$ over the category of $\mathsf{Q}_i$-schemes,
which to any $T\to \mathsf{Q}_i$ associates the set of
morphisms of schemes $f:\mathsf{G}_i\to G_i'\times_K\mathsf{Q}_i$ over
$\mathsf{Q}_i$, such that (i) the restriction to the identity component
of $\mathsf{G}_i$ is a homomorphism of group schemes, and
(ii) $f$ commutes with the right actions of identity components
of $\mathsf{G}_i$ and  $G_i'\times_K\mathsf{Q}_i$ on
$\mathsf{G}_i$ and  $G_i'\times_K\mathsf{Q}_i$ respectively.
Note that the identity component
of $\mathsf{G}_i$ is isomorphic to $H\times_K\mathsf{Q}_i$ as a group scheme.
In addition, since we have an isomorphism
$\mathsf{G}_i\cong G_i'\times_K\mathsf{Q}_i$ of schemes,
by choosing a closed point $c_j$ on each connected component
of $G_i'$ we can take sections $s_j:\mathsf{Q}_i\to \mathsf{G}_i$.
The functor
$\Hom^{\flat}_{\mathsf{Q}_i}(\mathsf{G}_i,G_i'\times_K\mathsf{Q}_i)$
is represented by a scheme over
$\mathsf{Q}_i$.
Indeed, a morphism $f:\mathsf{G}_i\to G_i'\times_K\mathsf{Q}_i$
with the properties (i) and (ii) is uniquely determined by
the restriction of $f$ to the connected component of $\mathsf{G}_i$
and the image of sections $\{s_j\}_j$ under $f$.
Therefore, if $\mathsf{G}_{i}^0$ denotes the identity component
of $\mathsf{G}_i$, then $\Hom_{\mathsf{Q}_i}(\mathsf{G}_{i}^0,G_i'\times_K\mathsf{Q}_i)\times_{\mathsf{Q}_i}(\mathsf{Q}_i\times_K(G_i')^r)$
represents $\Hom^{\flat}_{\mathsf{Q}_i}(\mathsf{G}_i,G_i'\times_K\mathsf{Q}_i)$, where $\Hom_{\mathsf{Q}_i}(\mathsf{G}_{i}^0,G_i'\times_K\mathsf{Q}_i)$
is the hom scheme which to any $T\to \mathsf{Q}_i$ associates
the set of group homomorphisms $\mathsf{G}_{i}^0\times_{\mathsf{Q}_i}T\to G_i'\times_KT$, and $r+1$ is the number of the connected components of $\mathsf{G}_i$.
By \cite[Expose XXIV 7.2.3]{SGA3}, $\Hom_{\mathsf{Q}_i}(\mathsf{G}_{i}^0,G_i'\times_K\mathsf{Q}_i)$ is a scheme (locally of finite type and separated) over
$\mathsf{Q}_i$. Thus, $\Hom^{\flat}_{\mathsf{Q}_i}(\mathsf{G}_i,G_i'\times_K\mathsf{Q}_i)$
is represented by a scheme over
$\mathsf{Q}_i$.
The morphism $\phi_i$ induces an inductive system of sections
$\{\Spec R_i/\mathfrak{m}_i^{n+1}\to \Hom^{\flat}_{\mathsf{Q}_i}(\mathsf{G}_i,G_i'\times_K\mathsf{Q}_i)\}_{n\ge 0}$, where $\Spec R_i=\mathsf{Q}_i$ and $\mathfrak{m}_i$ is the maximal ideal of $R_i$.
Since $\Hom^{\flat}_{\mathsf{Q}_i}(\mathsf{G}_i,G_i'\times_K\mathsf{Q}_i)$
is a scheme, the system is uniquely extended to a section
$\mathsf{Q}_i\to\Hom^{\flat}_{\mathsf{Q}_i}(\mathsf{G}_i,G_i'\times_K\mathsf{Q}_i)$. It gives rise to a morphism 
$\Phi_i:\mathsf{G}_i\to G'_i\times_K\mathsf{Q}_i$
that is an extension of $\phi_i$.
To check that $\Phi_i$ is a group homomorphism, it suffices to
show that $\Phi_i$ commutes with respect to group structures.
It follows from the facts that $\phi_i$ is a group homomorphism
and the natural completion map $\Gamma(\mathsf{G}_i\times_{\mathsf{Q}_i}\mathsf{G}_i,\OO_{\mathsf{G}_i\times_{\mathsf{Q}_i}\mathsf{G}_i})\to \Gamma(\mathsf{G}_i\times_{\mathsf{Q}_i}\mathsf{G}_i,\OO_{\mathsf{G}_i\times_{\mathsf{Q}_i}\mathsf{G}_i})^{\wedge}$ is injective (consider the compatibility of $\Phi_i$ with group structures
 in terms of ring homomorphisms).
For any $i$ and $j$, the intersection $(\mathsf{G}_i\times_{\mathsf{Q}_i}(\mathsf{Q}_i\cap \mathsf{Q}_j))\cap (\mathsf{G}_j\times_{\mathsf{Q}_j}(\mathsf{Q}_i\cap\mathsf{Q}_j))$ is a (reductive) smooth, affine group scheme over $\mathsf{Q}_i\cap\mathsf{Q}_j$. Two morphisms $\phi_i$ and $\phi_j$
coincide in the intersection (after associating the formal schemes),
and thus taking into account the above argument, we see
that $\Phi_i$ and $\Phi_j$ coincide in
$(\mathsf{G}_i\times_{\mathsf{Q}_i}(\mathsf{Q}_i\cap \mathsf{Q}_j))\cap (\mathsf{G}_j\times_{\mathsf{Q}_j}(\mathsf{Q}_i\cap\mathsf{Q}_j))$.
Therefore, $\Phi_i$'s are glued together, and thus
we have an isomorphism $\mathsf{G}\to H\times_K\mathsf{Q}$
of group schemes over $W^{\diamond}$.
\QED

(The second setup): Taking into account Lemma~\ref{nbdconst} and
Proposition~\ref{embed},
assume further that there exists a lifting $W\to \XX_U$ of $W\to \YY_U\cong [W/\bar{G}]$
such that:
\begin{enumerate}
\renewcommand{\labelenumi}{(\roman{enumi})}

\item $U$ is an affine scheme of finite type over $L$,

\item the pullback of $\mathcal{F}\to \XX$
by the composite $W\to \XX_U\to \XX$ is isomorphic to $W\times_{L}H$,

\item $(W\times_{\XX}I\XX)\times_{(W\times_L\bar{G})}(W\times_L\bar{G}')$ is embedded into $W\times_{L}G'$
as the inverse image of
$(W\times_{\YY_U}I\YY_U)\times_{(W\times_L\bar{G})}(W\times_L\bar{G}')\subset W\times_L\bar{G}'$ under $W\times_LG'\to W\times_L\bar{G}'$.
\end{enumerate}

\begin{Proposition}
\label{reducedquotgit}
With the notation as above, $\XX_U$ has the form
$[W/G]$.
\end{Proposition}

\Proof
Note first that
we have the lifting $\xi:W\to \XX_U$
and thus by Theorem~\ref{rigidification} (iii)
there exists a natural isomorphism
$\alpha:W\times_{\YY_U}\XX_U\to W\times_LBH$ over $W$.
This isomorphism is described as follows.
Let $\rho:\XX_U\to \YY_U$ be the rigidifying morphism.
For any $\omega:T\to W$, an object which belongs to $W\times_{\YY_U}\XX_U(T)$
amounts to data
\[
\{\omega:T\to W, \eta:T\to \XX_U,\ \theta:T\to \Isom_{\YY_U,T}(\rho(\eta),\bar{\xi}(\omega))\cong \Isom_{\XX_U,T}(\eta,\xi(\omega))/H\}
\]
where $\bar{\xi}$ is the composite $W\to \XX_U\to \YY_U$.
Then consider the principal $H$-bundle
\[
\Isom_{\XX_U,T}(\eta,\xi(\omega))\times_{\Isom_{\XX_U,T}(\eta,\xi(\omega))/H,\theta}T\to T
\]
and denote by $\phi_{\theta}:T\to BH$ the corresponding morphism.
Then $\alpha$ sends $(\omega,\eta,\theta)$ to $(\omega,\phi_{\theta})$.
There exists the diagram
\[
\xymatrix{
W\times_L\bar{G} \ar[r] \ar[d] & W\times_LBH \ar[r] \ar[d] & W \ar[d] \\
W \ar[r] & \XX_U \ar[r] & \YY_U \\
}
\]
where the right vertical arrow is determined by
$W\to [W/\bar{G}]$ and  all squares are cartesian.
Note that $W\times_L\bar{G}\to W\times_LBH$
is a morphism over $W$.
We will prove the claim:

\begin{Claim}
\label{factorclaim}
The composite $W\times_L\bar{G}\to W\times_LBH\to BH$
factors through  second projection
$W\times_L\bar{G}\to \bar{G}$.
\end{Claim}

Before the proof of the above claim, assuming that the claim holds,
we will complete the proof
of our Proposition.
By the claim, $W\times_L\bar{G}\to W\times_LBH$
arises from $\bar{G}\to BH$.
The fiber of $W\times_L\bar{G}\to W\times_LBH$
over $\Spec L\to W\times_L\Spec L\to W\times_L BH$
is $G\to \Spec L$. Here the first morphism $\Spec L\to W\times_L\Spec L$
is determined by
the unique point on $W$ lying over $u\in U$, and
the second morphism $W\times_L\Spec L\to W\times_L BH$ is defined
by the identity of $W$ and the natural projection $\Spec L\to BH$.
Thus we obtain the diagram
\[
\xymatrix{
W\times_LG \ar[r] \ar[d]^{\textup{pr}_1} & W\times_L\bar{G} \ar[r] \ar[d] & W \ar[d] \\
W \ar[r] & W\times_LBH \ar[r] & \XX_U \\
}
\]
where all squares are cartesian.
This means that $\XX_U$ is isomorphic to $[W/G]$
where the action of $G$ on $W$ is an extension of that of $\bar{G}$.
\QED

{\it Proof of Claim~\ref{factorclaim}}.
To prove our Claim, it suffices to prove that
the restriction of $W\times_{L}\bar{G}\to W\times_LBH\to BH$
to $W\times_{L}\bar{G}'$ factors through
the second projection $W\times_{L}\bar{G}'\to \bar{G}'$.
To see this reduction, it is enough to show that the
morphism $\phi:W\times_{L}\bar{G}'\to BH$
uniquely extends to a morphism $W\times_{L}\bar{G}\to BH$,
and any morphism $\psi:\bar{G}'\to BH$ uniquely extends to a
morphism $\bar{G}\to BH$.
To see this, 
set $\mathcal{I}$ and $\mathcal{J}$ are the ideals of $\OO_{W\times_{L}\bar{G}}$ and $\OO_{\bar{G}}$ determined by
$W\times_{L}\bar{G}'$ and $\bar{G}'$ respectively,
and assume that $\mathcal{I}^2=0$ and $\mathcal{J}^2=0$.
By the deformation theory \cite[Theorem 1.5]{OL},
it suffices to check $\Ext^0(L\phi^*\mathsf{L}_{BH/\Spec L},\mathcal{I})=\Ext^1(L\phi^*\mathsf{L}_{BH/\Spec L},\mathcal{I})=0$ and
$\Ext^0(L\psi^*\mathsf{L}_{BH/\Spec L},\mathcal{J})=\Ext^1(L\psi^*\mathsf{L}_{BH/\Spec L},\mathcal{J})=0$.
Notice that $\phi$ and $\psi$ are flat.
(Any morphism $X\to BH$ factors through the natural
projection $\Spec L\to BH$ after replacing $X$ with some fppf cover $X'\to X$.)
Thus by the same argument in the proof of Lemma~\ref{liftgit},
it suffices to see that $\mathsf{L}_{BH/\Spec L}$ is of perfect amplitude
in $[1]$. It follows from Remark~\ref{smoothcotangent} (i).

Next observe the morphism $W\times_L\bar{G}\to W\times_{\YY_U}\XX_U\cong W\times_LBH$.
Let $(\omega,g)$ be a $T$-valued point of $W\times_L\bar{G}$,
where $\omega:T\to W$ and $g:T\to \bar{G}$.
Considering the cartesian diagram
\[
\xymatrix{
W\times_L\bar{G} \ar[r] \ar[d] & W \ar[d] \\
W\times_{\YY_U}\XX_U \ar[r] & \XX_U \\
}
\]
where the top horizontal arrow is the action of $\bar{G}$,
we see that the image of $(\omega,g)$
is described by
\[
\{ \omega\in W(T), g\omega\in \XX_U(T), g\in \Isom_{\YY_U,T}(\bar{\omega},\bar{g\omega})(T)\}
\]
where $\bar{\omega}$ and
$\bar{g\omega}$ denote the images of $\omega$ and $g\omega$
in $\YY_U$ respectively, and for any $f:T'\to T$ we identify $\Isom_{\YY_U,T}(\bar{\omega},\bar{g\omega})(T')$
with $\{g'\in \bar{G}(T')|\ g'(f^*\omega)=f^*(g\omega)\}$.
Thus we have the diagram
\[
\xymatrix{
\mathcal{P}\ar[r] \ar[d] & \Isom_{\XX_U,T}(\omega,g\omega) \ar[d] & \\
T \ar[r] & \Isom_{\YY_U,T}(\bar{\omega},\bar{g\omega})\ar[r] & \bar{G}\times_LT \\
}
\]
where the square is cartesian, and the lower right horizontal arrow
is a closed immersion.
The principal $H$-bundle $\mathcal{P}$ corresponds to
the image of $W\times_{L}\bar{G}\to W\times_LBH\to BH$.
Suppose $g\in \bar{G}'(T)$.
Let $\Isom_{\YY_U,T}'(\bar{\omega},\bar{g\omega}):=\Isom_{\YY_U,T}(\bar{\omega},\bar{g\omega})\times_{\bar{G}\times_LT}(\bar{G}'\times_LT)$
and $\Isom_{\XX_U,T}'(\omega,g\omega):=\Isom_{\XX_U,T}(\omega,g\omega)\times_{\Isom_{\YY_U,T}(\bar{\omega},\bar{g\omega})}\Isom_{\YY_U,T}'(\bar{\omega},\bar{g\omega})$.
Here we claim that $\Isom_{\XX_U,T}'(\omega,g\omega)\to \Isom_{\YY_U,T}'(\bar{\omega},\bar{g\omega})$ can be identified
with the principal $H$-bundle
$\Isom_{\YY_U,T}'(\bar{\omega},\bar{g\omega})\times_{(\bar{G}'\times_LT)}(G'\times_LT)
\to \Isom_{\YY_U,T}'(\bar{\omega},\bar{g\omega})$.
In particular, $\mathcal{P}$ depends only on $g\in \bar{G}'(T)$.
To this end, we first consider the case when
$g:T\to \bar{G}'$ extends to $\tilde{g}:T\to G'$ and $\omega\cong g\omega$ in $\XX_U$.
If $\Stab'(\omega):=\Isom_{\YY_U,T}'(\bar{\omega},\bar{\omega})\subset \bar{G}'\times_LT$,
then $\Isom_{\YY_U,T}(\bar{\omega},\bar{g\omega})$
is $g\cdot \Stab'(\omega)$ in $\bar{G}'\times_LT$.
Thus if $G'$ acts on $W$ via $G'\to \bar{G}'$
and $\widetilde{\Stab'(\omega)}$ denotes the stabilizer group scheme
of $\omega$ with respect to the action of $G'$, then
the inverse image of $\Isom_{\YY_U,T}'(\bar{\omega},\bar{g\omega})$
in $G'\times_LT$ is $\tilde{g}\cdot \widetilde{\Stab'(\omega)}$.
By our assumption and the second setup (iii), we have $\widetilde{\Stab'(\omega)}\cong \Aut_{\XX_U,T}(\omega)\times_{\Isom_{\YY_U,T}(\bar{\omega},\bar{\omega})}\Isom_{\YY_U,T}'(\bar{\omega},\bar{\omega})$
over $T$.
In the general case,
take an \'etale surjective morphism $T'\to T$
so that $T'\to T\to \bar{G}'$ extends to $\tilde{g}:T'\to G'$ and $w|_{T'}\cong gw|_{T'}$ in $\XX_{U}$.
Then by the previous case and the descent theory,
there exists a dotted closed immersion
\[
\xymatrix{
\Isom_{\XX_U,T'\times_TT'}'(\omega|_{T'\times_TT'},\tilde{g}\omega|_{T'\times_TT'}) \ar@<0.5ex>[r]^(0.6){\textup{pr}_1} \ar@<-0.5ex>[r]_(0.6){\textup{pr}_2}  \ar[d] & \Isom_{\XX_U,T'}'(\omega|_{T'},\tilde{g}\omega|_{T'})\ar[r] \ar[d] & \Isom_{\XX_U,T}'(\omega,\tilde{g}\omega) \ar@{..>}[d]\\
G'\times_L(T'\times_TT') \ar@<0.5ex>[r]^{\textup{pr}_1} \ar@<-0.5ex>[r]_{\textup{pr}_2}  & G'\times_LT' \ar[r] & G'\times_LT, \\
}
\]
where the left and central vertical arrows are pullbacks of 
\[
\Isom_{\YY_U,T'\times_TT'}'(\bar{\omega}|_{T'\times_TT'},\bar{g\omega}|_{T'\times_TT'})\subset \bar{G}'\times_L(T'\times_TT')
\]
and $\Isom_{\YY_U,T'}'(\bar{\omega}|_{T'},\bar{g\omega}|_{T'})\subset \bar{G}'\times_LT'$ respectively.
This implies our Claim.
\QED

\begin{Theorem}
\label{maingit}
Let $\XX$  be an algebraic stack locally of finite type over a perfect field.
Then the open substack $\XX^{\gs}$ of GIT-like p-stable points
has a coarse moduli map $\pi:\XX^{\gs}\to X$,
such that $X$ is locally of finite type.
Moreover, $\pi$ is universally closed and quasi-finite,
which induces an isomorphism $\OO_X\to \pi_*\OO_{\XX}$.
If $X'\to X$ is a flat morphism, then the
second projection $\XX^{\gs}\times_{X}X'\to X'$ is a coarse moduli map.
\end{Theorem}

\Proof
Let $\XX_0^{\gs}$ be the reduced algebraic stack
associated to $\XX^{\gs}$.
By Proposition~\ref{reducedgit}
there exists a coarse moduli map $\pi_0:\XX_0^{\gs}\to X_0$.
Then by Proposition~\ref{reducedquotgit},
any point on $X_0$ has an \'etale neighborhood
$U_0\to X_0$, such that $\XX_0^{\gs}\times_{X_0}U_0$ has
the form $[W_0/G]$, where $W_0$ is an affine scheme which is finite over
the affine scheme $U_0$, and $G$ is a linearly reductive group.
By Lemma~\ref{liftgit}, the \'etale deformation $\XX_U\to \XX^{\gs}$
of
$\XX_0^{\gs}\times_{X_0}U_0\to \XX_0^{\gs}$ to $\XX^{\gs}$ has the form $[W/G]$.
Since $G$ is linearly reductive,
thus
\[
\Gamma(\XX_U,\OO_{\XX_U})=\Gamma(W,\OO_W)^G\to \Gamma(\XX_0^{\gs}\times_{X_0}U_0,\OO_{\XX_0^{\gs}\times_{X_0}U_0})=\Gamma(W_0,\OO_{W_0})^G
\]
is surjective.
Applying Lemma~\ref{helpgit}
we conclude that the nilpotent deformation
$\XX_0^{\gs}\to \XX^{\gs}$ has the property $(\mathbb{L})$ in section 4.
Therefore according to Proposition~\ref{deformmoduli}$,
\XX^{\gs}$ has a coarse moduli map.
By the construction, it is clear that $\pi$ is
universally closed and quasi-finite.
The last assertion also follows from Proposition~\ref{deformmoduli}.
\QED

From Proposition~\ref{deformmoduli}, Theorem~\ref{maingit}, Lemma~\ref{liftgit}, Proposition~\ref{reducedquot}, we deduce

\begin{Corollary}
\label{gitlocal}
Let $\XX$ be an algebraic stack locally of finite type over a perfect field.
Let $\pi:\XX^{\gs}\to X$ be a coarse moduli map.
Then every point on $X$ admits an \'etale neighborhood
$U\to X$ such that
$\XX^{\gs}\times_{X}U\to U$
has the form $[W/G]$ where $W$ is finite
over $U$, and $G$ is a linearly reductive group acting on $W$ over $U$.
More precisely, $U$ is an affine scheme over a finite (separable) extention $k'\supset k$ of fields
and $G$ is a linearly reductive group over $k'$, which acts
on $W$ over $U$ so that the quotient stack $[W/G]$ is isomorphic to
$\XX^{\gs}\times_{X}U$ over $U$.
\end{Corollary}

\begin{Remark}[Isovariant \'etale]
Corollary~\ref{gitlocal} says that
GIT-like stable point is approximated by a quotient stack via
an \'etale morphism $[W/G]\to \XX$ (with the above notation).
Moreover, $[W/G]\to \XX$ preserves the structures of automorphisms.
More precisely, if $I[W/G]\to [W/G]$ and $I\XX\to \XX$ denote
the inertia stacks of $[W/G]$ and $\XX$ respectively,
then the natural morphism $I[W/G]\to I\XX\times_{\XX}[W/G]$
is an isomorphism.
It is quite useful.
For instance, generalizing
Thomason's descent theory (\cite{Thomason}) Joshua developed the
{\it isovariant \'etale
descent theory} for G-theory on algebraic stacks
(\cite[Section 5]{Joshua}).
Informally, an isovariant \'etale morphism of algebraic stacks is
an \'etale morphism which preserves the structures of automorphisms (cf. \cite[DEFINITIONS 3.1 (iii)]{Joshua}), and the notion of
isovariant \'etaleness is crucial for
the descent theory of G-theory.
Note that the above morphism $[W/G]\to \XX$ is isovariant \'etale.
Moreover, under the stable condition (see Definition~\ref{GIT-likestabletype}
, Remark~\ref{zerodef}),
one can reduce G-theory of algebraic stacks to equivariant G-theory, that is,
the case of an affine scheme provided with the action of a reductive group.
Since this topic is beyond the scope of this paper,
we will discuss these issues and applications in another paper.
\end{Remark}

Motivated by Theorem~\ref{maingit} and Corollary~\ref{gitlocal},
we propose a class of Artin stacks.

\begin{Definition}
\label{GIT-likestabletype}
Let $\XX$ be an algebraic stack locally of finite type
over a perfect field $k$.
We say that $\XX$ is {\it of GIT-like stable type} over $k$
(or simply {\it stable algebraic stack} over $k$)
 if $\XX^{\gs}=\XX$.
In other words, $\XX$ is of GIT-like stable type
if the automorphism groups of all closed points are linearly reductive
and any closed point admits an open neighborhood $\UU\subset \XX$,
such that the inertia stack $I\UU_0$ of the reduced stack associated to $\UU$
has a closed subgroup $\FF_0$
that satisfies:
\begin{enumerate}
\renewcommand{\labelenumi}{(\roman{enumi})}

\item $\FF_0$ is smooth and affine over $\UU_0$, and all geometric fibers are
connected,

\item the quotient $I\UU_0/\FF_0$ is finite.

\end{enumerate}

\end{Definition}

\begin{Remark}
\label{zerodef}
By Remark~\ref{gitdefremark}, in the case of
characteristic zero $\XX$ is of GIT-like stable type if and only if
the followings hold: (i) the automorphism group of every geometric
point on $\XX$ is reductive,
(ii) if $\XX_{\operatorname{red}}$ denotes the reduced stack
associated to $\XX$ then the
inertia stack $I\XX_{\operatorname{red}}$ is Zariski locally equidimensional over $\XX_{\operatorname{red}}$,
(iii) the identity component of $I\XX_{\operatorname{red}}$ is affine over $\XX_{\operatorname{red}}$, and
(iv) the quotient of $I\XX_{\operatorname{red}}$ by its identity component
is finite over $\XX_{\operatorname{red}}$. (See Theorem~\ref{characterization}.)
\end{Remark}

In virtue of the works of Inaba \cite{Inaba}
Lieblich \cite{lieb} and To\"en-Vaqui\'e \cite{TV},
we have
the moduli Artin stack $\mathfrak{D}^{b}_{\textup{p}}(X)$ of objects in the derived category of perfect
complexes (satisfying a certain condition) on a proper flat scheme $X$.
Also, the general result of \cite{TV} yields the moduli stack
of complexes of representations of a finite quiver.
Recent developments on derived category reveal the importance of
these moduli stacks.
On the other hand, we would like to call attention to the fact:
one cannot interpret these stacks 
as quotient stacks (at least a priori since abstract approaches are applied).
It would be interesting to apply our GIT-like p-stability to these stacks.
We hope to come back to this topic in a future work.

%
%
%
%
%
%
%
%

\section{Comparing with Mumford's Geometric Invariant Theory}
By an algebraic scheme over a field $k$ we mean a scheme 
locally of finite type and separated over a field $k$.
In this section, we assume that the base field $k$
is algebraically closed of characteristic zero except Theorem~\ref{characterization}.
In this section we discuss the relationship between
our GIT-like p-stability
and Geometric Invariant Theory due to Mumford (\cite{GIT}).
Moreover, we prove Theorem~\ref{B}.

Let $X$ be an algebraic scheme over $k$.
Let $G$ be a linearly reductive group scheme over $k$.
(An algebraic group over $k$ is linearly reductive if and only if
it is reductive.)
Let $\sigma:G\times_kX\to X$ be an action on $X$.
Let $X(\textup{Pre})\subset X$ be the open subset of $X$, consisting of pre-stable points
in the sense of \cite[Definition 1.7]{GIT}.
The main purpose of this section is to prove:

\begin{Theorem}
\label{gitrelation}
Let $[\XPre/G]$ be the open substack of $[X/G]$,
associated to $\XPre$.
Let $[X/G]^{\gs}$ be the open substack consisting of GIT-like
p-stable points on $[X/G]$.
Let $\mathcal{S}$ be the maximal open substack of $[X/G]$, 
admitting a coarse moduli space that is a scheme.
Then 
\[
[\XPre/G]=[X/G]^{\gs}\cap \mathcal{S}.
\]
\end{Theorem}

\begin{Remark}
\label{Crem}
We would like to invite your attention to some of the advantages
of our approach.
First of all, as in Keel-Mori theorem coarse moduli are allowed
to be algebraic spaces (rather than schemes).
It makes the framework more amenable.
Moreover, contrary to Geometric Invariant Theory
our approach does not rely on global quotient structures (cf. Remark~\ref{zerodef}).
It is important for several reasons.
First, in practice it is hard to prove that a give algebraic stack
is a quotient stack.
In addition, algebraic stacks do not necessarily have such structures
(see \cite[section 2]{EHKV}).
Secondly, as mentioned at the close of section 5,
nowadays we often use abstract methods for constructing algebraic moduli
stacks, such as Artin's representability theorem, and
Geometric Invariant Theory is not applicable to stacks constructed by such abstract methods. For example, Lieblich applied the Artin's theorem to the
constructions of moduli stacks of twisted sheaves (\cite{lieb2}) and complexes
(\cite{lieb}), and Olsson used the theorem in the work on Hom stacks (\cite{Hom}). Also, Lurie proved an amazing generalization of Artin's representability theorem
to derived algebraic geometry (\cite{Lu}).
\end{Remark}

\begin{Remark}
According to Theorem~\ref{gitrelation}, we may say that in a sense 
the notion of GIT-like p-stability is an intrinsic generalization
of pre-stability, that is, the ``local part'' of Geometric
Invariant Theory.

By the theory of good moduli spaces \cite[Theorem 6.6]{Alper}, the geometric quotient of $\XPre$ by $G$ in the sense of Definition 0.6 of \cite{GIT}
is a coarse moduli space for $[\XPre/G]$.
(We will recall the notion of good moduli spaces introduced
by Alper in section 7.)
Thus $[\XPre/G]\subset \mathcal{S}$.
Namely, $[\XPre/G]=[\XPre/G]\cap \mathcal{S}$.
\end{Remark}

We first show $[\XPre/G]\subset [X/G]^{\gs}$. Note that if $[\XPre/G]$
is contained in $[X/G]^{\gs}$, then $[\XPre/G]\subset [X/G]^{\gs}\cap \mathcal{S}$.

\begin{Proposition}
\label{relation1}
Any closed point on $[\XPre/G]$
is a GIT-like p-stable point.
\end{Proposition}

To prove this Proposition,
we may and will assume that $X$ is affine,
and the action of $G$ on $X$ is closed (cf.
\cite[Definition 1.7]{GIT}).

\begin{Lemma}
\label{matsushima}
Let $p$ be a closed point on $\XPre$.
Then the stabilizer group is (linearly) reductive.
\end{Lemma}

\Proof
According to Matsushima's theorem (cf. \cite[page 84]{L}),
an orbit is affine if and only if the stabilizer is reductive.
Since every orbit in $\XPre$ is a closed set in an affine open set,
thus our claim follows.
\QED

Let $X_{0}$ be the reduced affine scheme associated to $X$.
The base field is perfect, and thus
$G\times_kX_0$ is also reduced.
Therefore, the action on $X$ induces an (closed) action
$\sigma_0:G\times_kX_0\to X_0$.
Let $\mathsf{Stab}\to X_0$
be the stabilizer group, which is defined to be the second projection
$(G\times_kX_0)\times_{X_0\times_kX_0}X_0\to X_0$,
where $(\sigma_0,\textup{pr}_{2}):G\times_kX_0\to X_0\times_k X_0$,
and $X_0\to X_0\times_k X_0$ is the diagonal map.
Let $\mathcal{F}$ be the identity component of
$\mathsf{Stab}$.
Note that $\mathcal{F}\to X_0$ is Zariski locally equidimensional (cf. \cite[page 10]{GIT})
and geometric fibers are connected and smooth because the base field is characteristic zero.
Then by \cite[$VI_B$ Corollarie 4.4]{SGA3}, 
$\mathcal{F}\to X_0$ is a smooth open normal subgroup of $\mathsf{Stab}$.
Note that the inertia stack $I[X_0/G]\to [X_0/G]$
is described by
\[
[\mathsf{Stab}/G]\to [X_0/G],
\]
where $G$ acts on $\mathsf{Stab}\subset G\times_kX_0$ by $(h,x)\mapsto (ghg^{-1},gx)$
for any $(h,x)\in \mathsf{Stab}\subset G\times_kX_0$ and any $g\in G$.
Thus $\mathcal{F}$ descends to an open subgroup
stack $[\mathcal{F}/G]\to [X_0/G]$ of the inertia stack.
Moreover, the affineness of $\mathsf{Stab}$ over $X_0$ together with
the following Lemma~\ref{finitefiltration} implies that $\mathcal{F}$
is affine over $X_0$.

\begin{Lemma}
\label{finitefiltration}
The quotient $\mathsf{Stab}/\mathcal{F}$ is finite
over $X_0$.
\end{Lemma}

\Proof
Let $x\in X_0$ be a closed point on $X$.
Let $G_x$ be the stabilizer group of $x$.
The \'etale slice theorem of Luna (\cite[page 96]{OL}) says
that there exist a locally closed affine $G_x$-invariant
subscheme $V\subset X_0$ containing $x$,
and the commutative diagram
\[
\xymatrix{
(V\times G)/G_x \ar[r]  \ar[d] & [V/G_x] \ar[r] \ar[d] & Z \ar[d] \\
X_0 \ar[r]  & [X_0/G] \ar[r] & Y, \\
}
\]
where all squares are cartesian, all vertical arrows are affine
and \'etale, and $Z$ (resp. $Y$) are geometric quotients of $V$
(resp. $X_0$) by $G_x$ (resp. $G$).
Here if $V=\Spec A$ and $X_0=\Spec B$, then
$Z=\Spec A^{G_{x}}$ and $Y=\Spec B^{G}$,
and $(V\times_kG)/G_x$ is the quotient of $V\times_kG$
by the action of $G_x$ determined by $h\cdot (v,g)= (hv,gh^{-1})$.
Using the diagram, we want to reduce the problem to $[V/G_x]$.
Let $\mathsf{Stab}_V\to V$ be the stabilizer group scheme of the action
of $G_x$ on $V$. Since $\mathsf{Stab}_V$ is equidimensional, thus
$\mathsf{Stab}_V$ contains $G_{x}^0\times_kV$, where $G_{x}^0$
is the identity component of $G_x$.
To see that $\mathsf{Stab}/\mathcal{F}$ is finite over $X_0$,
it suffices to show that $\mathsf{Stab}_V/(G_{x}^0\times_kV)$
is finite over $V$.
We can consider the action of $G_{x}$ on $V$ to be the action
of the finite group $\bar{G}_x:=G_{x}/G_{x}^0$.
Then it follows from \cite[Proposition 0.8]{GIT}
that the action of $\bar{G}_x$ on $V$ is proper.
In particular, $\mathsf{Stab}_V/(G_{x}^0\times_kV)$
is finite over $V$. This completes the proof.
\QED

{\it Proof of Proposition~\ref{relation1}.}
By Lemma~\ref{matsushima},
every closed point on $[\XPre/G]$ has a (linearly) reductive automorphism
group.
Let $x\in [X/G]=[\XPre/G]$ be a closed point.
Take a closed point $x'\in X$ lying over $x$.
Let $\hat{\OO}_{X,x'}$ be the completion of local ring $\OO_{X,x'}$ at $x'$.
Then $\Spec \hat{\OO}_{X,x'}\to [X/G]$ is an effective versal deformation for $x$.
And by Lemma~\ref{finitefiltration}, the restriction of the subgroup scheme $\mathcal{F}$ to $\Spec \hat{\OO}_{X,x'}$ yields
(b) in
Definition~\ref{defgit}.
\QED

Next we will show the converse:

\begin{Proposition}
\label{relation2}
We have $[\XPre/G]\supset [X/G]^{\gs}\cap \mathcal{S}$.
\end{Proposition}

\Proof
We say that a point $x\in X$ is regular if it admits an open
 neighborhood
on which stabilizers have constant dimension.
Let $X^{\textup{reg}}$ denote the open set of regular points.
Clearly, $[X/G]^{\gs}\subset [X^{\textup{reg}}/G]$.
Let $V\subset X$ be the open subset consisting of points lying over $[X/G]^{\gs}$.
Let $v:V\to  [X/G]^{\gs}$ be the natural projection.
Let $p:[X/G]^{\gs}\to Z$ be a coarse moduli map.
Let $z'\in [X/G]^{\gs}$ be a closed point and $z$ the image in $Z$.
Suppose that $z$ has a Zariski open neighborhood
which is a scheme.
We will show that $z$ has an affine neighborhood $T$ such that
$v^{-1}p^{-1}(T)$ is an affine scheme.
By Corollary~\ref{gitlocal}, there exists
an \'etale neighborhood $U\to Z$ of $z$,
such that $[X/G]^{\gs}\times_{Z}U\to U$ has the form $[W/H]\to U$
where $W$ is affine over an affine scheme $U$, and $H$ is
a linearly reductive group over $k$.
Then $[W/H]\to BH\times_k U$ is affine, and $BH\times_kU\to U$
gives rise to an exact push forward functor from the category
of quasi-coherent sheaves
on $BH\times_kU$ to that of quasi-coherent sheaves
on $U$.
If $T'$ denotes the image of $U\to Z$, then we can conclude that
$(p\circ v)_*$ is an exact functor from the category of quasi-coherent
sheaves on $v^{-1}p^{-1}(T')$ to that of quasi-coherent sheaves on $T'$.
Thus $v^{-1}p^{-1}(T')\to p^{-1}(T')\to T'$ is an affine morphism
(notice that $v:V\to [X/G]^{\gs}$ is a principal $G$-bundle, in particular
an affine morphism)
. Take an affine neighborhood $z\in T\subset T'$.
Then $v^{-1}p^{-1}(T)$ is affine.
Since $[X/G]^{\gs}\subset [X^{\textup{reg}}/G]$,
the action on $v^{-1}p^{-1}(T)$
is closed by \cite[page 10, line 19-20]{GIT}.
This implies our claim.
\QED

{\it Proof of Theorem~\ref{gitrelation}}.
Proposition~\ref{relation1} and~\ref{relation2} imply
our claim.
\QED

\begin{Remark}
In the proof of Proposition~\ref{relation1},
we can find the filtration
as in Definition~\ref{defgit} (b).
\end{Remark}

In the proof of Lemma~\ref{matsushima} and~\ref{finitefiltration},
we observed the closed action of a (linearly) reductive group
on an affine scheme.
From the observation, we have the following characterization:

\begin{Theorem}
\label{characterization}
Let $\XX$ be an algebraic stack locally of finite type
over a field $k$ of characteristic zero.
Then $\XX$ is stable algebraic stack (Definition~\ref{GIT-likestabletype} and Remark~\ref{zerodef}) if and only if
the following conditions hold:
\begin{enumerate}
\renewcommand{\labelenumi}{(\roman{enumi})}
\item There exists a coarse moduli map $\pi:\XX\to X$,
such that $X$ is locally of finite type.
For any \'etale morphism $X'\to X$ of algebraic spaces,
the second projection $\XX\times_{X}X'\to X'$ is a coarse moduli map.

\item For any point $x\in X$, there exists an \'etale neighborhood
$U\to X$, a finite extension of $k'\supset k$ fields, and a (linearly)
reductive group $G$ over $k'$,
such that $U$ is an affine  $k'$-scheme and
$\XX\times_{X}U$ has the form
$[W/G]$, where $W$ is a scheme that is affine over
$U$ and $G$ acts on $W$ over $U$.
\end{enumerate}
\end{Theorem}

\Proof
The ``only if'' direction follows from Theorem~\ref{maingit}
and Corollary~\ref{gitlocal}.
We next prove the ``if'' direction.
Taking into account Remark~\ref{gitdefremark} (iii)
we may and will assume that the base field is algebraically closed.
As noted above, by Lemma~\ref{matsushima} and~\ref{finitefiltration}
it suffices to show that if $[W/G]\to U$ is a quotient stack as in
(ii), then the action of $G$ is closed.
Since $[W/G]\to U$ is a coarse moduli map, the action of $G$ on $W$
is closed, that is, every orbit is closed. This completes the proof.
\QED

\begin{Remark}
Perhaps one wish to have a necessary and sufficient condition
for the existence of a coarse moduli space.
However, if one removes the condition (ii) Theorem~\ref{characterization},
then there are pathological examples even
in the case of Deligne-Mumford stacks (that have only quasi-finite inertia stacks).
For instance, it happens that a Deligne-Mumford stack
which does not have finite inertia stack, has a coarse moduli space.
Let $\mathbb{A}^1_{\CC}$ be a complex affine line and  $f:\mathbb{A}^1_{\CC}\sqcup \mathbb{A}^1_{\CC}\to \mathbb{A}^1_{\CC}$ the fold map.
Removing one point $q$ from $f^{-1}(0)$, we obtain
 a flat group scheme $g=f|_G:G=\{\mathbb{A}^1_{\CC}\sqcup \mathbb{A}^1_{\CC}\setminus q\}\to\mathbb{A}_{\CC}^1$ over $\mathbb{A}_{\CC}^{1}$, such that
 $g^{-1}(p)$ consists of two points (resp. one point) if $p\neq 0$
 (resp. $p=0$).
Let $BG$ be the classifying stack of $G$ over $\mathbb{A}^1$.
Then the inertia stack is quasi-finite over $BG$ but not finite over $BG$.
However, $\mathbb{A}^1_{\CC}$ is a coarse moduli space for $BG$.
We can not apply Keel-Mori theorem  to $BG$, and
in particular $BG^{\gs}$ does not coincide with $BG$.
The point is that $BG$ does not satisfies (ii) in
Theorem~\ref{characterization}.

Furthermore there exists a Deligne-Mumford stack
which does not admit a coarse moduli space.
Such an example can be found in Example 7.15 of \cite{Rydh}.
\end{Remark}

\begin{Remark}
Strong p-stability can differ from
Mumford's theory.
We will present the simple example which illustrates it.
Let $\mathbb{P}^1_{\CC}$ be the projective sphere over the
complex number field $\CC$.
Let $G=\operatorname{PGL}(2,\CC)$ be the algebraic group
that is the automorphism group of $\mathbb{P}_{\CC}^1$ over $\CC$.
The algebraic group $G$ acts on
$\mathbb{P}^1_{\CC}$ in the natural manner.
Every affine open set $U$ on $\mathbb{P}^1_{\CC}$ is not
$G$-invariant, that is, $G(U)\neq U$.
Thus every point on $\mathbb{P}^1_{\CC}$
is not pre-stable in the sense of Geometric Invariant Theory.
In fact, the stabilizer group is unipotent.
Let
$\mathsf{Stab}\to \mathbb{P}^1_{\CC}$ denote
the stabilizer group scheme over $\mathbb{P}^1_{\CC}$.
Then we can easily see that $\mathsf{Stab}\to \mathbb{P}^1_{\CC}$
is flat.
Therefore, every (and only one) closed point
on $[\mathbb{P}^1_{\CC}/G]$ is strong p-stable in our sense.
In this case, the coarse moduli space is $\Spec \CC$.
In general, algebraic stacks which have non-reductive 
(and positive dimensional) automorphisms are complicated.
For instance, by Nagata's example the finite generation of invariant rings
does not hold.

\end{Remark}

\section{Application}
In this section, we discuss the finiteness of coherent cohomology
as an application.

\begin{Proposition}
\label{behavelikeproper}
Let $\XX$ be a stable algebraic stack locally of finite type 
over a perfect field.
Let $\pi:\XX\to X$ be a coarse moduli map.
Let $\mathcal{F}$ be a coherent sheaf on $\XX$.
Then $R\pi_*^0\mathcal{F}$ is a coherent sheaf on $\XX$,
and $R\pi_*^i\mathcal{F}=0$ for $i>0$.
\end{Proposition}

\Proof
Our claim is \'etale local on $X$.
Thus by Corollary~\ref{gitlocal}
we may assume that $\XX\to X$ is of the form
$[W/G]\to U$, where $W=\Spec B$ is finite over an affine scheme $U=\Spec A$,
and
$G$ is a linearly reductive group.
A coherent sheaf $\mathcal{F}$ on $[W/G]$
amounts to a finite $B$-module $M$ which is equipped
with a $G$-action commuting with the action on $B$.
The direct image $\pi_*\mathcal{F}$ corresponds to $M^{G}$.
Thus the rule $\mathcal{F}\mapsto \pi_*\mathcal{F}$
is an exact functor because $G$ is linearly reductive.
Hence $R\pi_*^i\mathcal{F}=0$ for $i>0$.
It remains to prove that $R\pi_*^0\mathcal{F}$ is a coherent sheaf on $\XX$.
It is enough to show that $M^{G}$ is a finite $A$-module.
Note that $M$ is a finite $A$-module,
and $M^{G}$ is a sub $A$-module.
Since $A$ is noetherian, $M^{G}$ is a finite $A$-module.
\QED

\begin{Theorem}
\label{coherency}
Let $\XX$ be an algebraic stack locally of finite type
over a field $k$. Let $\bar{k}\supset k$ be an algebraic closure,
and suppose that $\XX\times_k\bar{k}$ is 
a stable algebraic stack over $\bar{k}$
and its coarse moduli space $X$ is proper over $\bar{k}$.
Let $\mathcal{F}$ be a coherent sheaf on $\XX$.
Then the cohomology group
$H^i(\XX,\mathcal{F})$ is finite dimensional for $i\ge 0$.
\end{Theorem}

\Proof
We may assume that the base field is algebraically closed.
Consider the composition $\XX\stackrel{\pi}{\to} X \to \Spec k$,
where $\pi$ is a coarse moduli map.
Taking into account the finiteness theorem for algebraic space
\cite[IV Theorem 4.1]{Kn} and
Leray spectral sequence for the composition,
it suffices to prove that
$R^i\pi_*\mathcal{F}$ is coherent for all $i$.
Thus our claim follows from Proposition~\ref{behavelikeproper}.
\QED

For example, it has the following direct corollary.

\begin{Corollary}
Let $\XX$ be an algebraic stack of finite type over a field $k$.
Suppose that $\XX\times_{k}\bar{k}$ is a stable algebraic stack,
where $\bar{k}$ is an algebraic closure of $k$.
If a coarse moduli space for $\XX\times_{k}\bar{k}$
is proper over $\bar{k}$, then $\XX$ has a versal deformation
(cf. \cite[Definition 4.1.1]{Aoki}).
\end{Corollary}

\Proof
By the proof of \cite[section 4.2]{Aoki}, it suffices to
show the finiteness of coherent cohomology.
Thus, our claim follows from Theorem~\ref{coherency}.
\QED

\begin{Remark}
If $\XX$ is proper over $k$, the above corollary is due to Aoki (cf. \cite[Theorem 1.3]{Aoki}).
\end{Remark}

\begin{Remark}
The finiteness theorem of coherent cohomology for proper algebraic stacks
has been proved by Laumon and Moret-Bailly under some hypotheses
(cf. \cite[(15.6)]{LM}).
Later,
Faltings proved the finiteness theorem for general proper stacks
via a surprising method of rigid geometry (cf. \cite{Fal}).
Recently, Olsson-Gabber proved Chow's lemma for algebraic stacks
and reproved the finiteness theorem (cf. \cite{OLp}).

We would like to stress that our theorem can be applied to
algebraic stacks whose stabilizers are linearly reductive and positive
dimensional.
To compare our result with the previous theorems,
consider a proper algebraic stack
$\XX\to \Spec A$.
Then the diagonal
$\XX\to \XX\times_{\Spec A}\XX$ is proper.
Practically, in many cases,
the proper diagonal is a {\it finite} morphism.
(In characteristic zero, it also implies that $\XX$ is Deligne-Mumford.)
Then if $\XX$ has finite diagonal, then by Keel-Mori theorem
$\XX$ has a proper coarse moduli space
$X$.
Therefore we may summarize the above as follows:
The finiteness theorem for proper stacks practically
tells us that if an algebraic stack has finite diagonal
and a proper coarse moduli
space, then it has the finiteness of coherent cohomology.
The main advantage of our finiteness result is
that it is applicable to a certain class of algebraic stacks whose stabilizers
are positive dimensional.
\end{Remark}

\begin{Remark}
Here we would like to relate our results with the notion of
good moduli spaces introduced by Alper (\cite{Alper}).
Let us recall the definition of good moduli space:
Let $\operatorname{QCoh}(\XX)$ denote the abelian category of
quasi-coherent sheaves on an algebraic stack $\XX$.
Let $f:\XX\to X$ be a morphism to an algebraic space $X$.
The morphism $f:\XX\to X$ is said to be a good moduli space for $\XX$
if the followings hold:
\begin{enumerate}
\renewcommand{\labelenumi}{(\alph{enumi})}

\item $f$ is a quasi-compact morphism and the push-forword functor
\[
f_*:\operatorname{QCoh}(\XX)\longrightarrow \operatorname{QCoh}(X)
\]
is exact,

\item the natural morphism $\OO_X\to f_*\OO_{\XX}$ is an isomorphism.
\end{enumerate}
In \cite{Alper}, many properties of good moduli spaces are
systematically studied.
It is worth remarking that by Theorem~\ref{maingit} and Proposition~\ref{behavelikeproper} a stable algebraic stack admits a good moduli space.
(Notice that the proof of Proposition~\ref{behavelikeproper} shows
the condition (a).)
\end{Remark}

Let $\XX$ be an algebraic stack of finite type over a perfect field
$k$. Suppose that $\XX$ is of GIT-like 
stable type over $k$.
Let $\pi:\XX\to X$ be a coarse moduli map (cf. Theorem~\ref{maingit}).
The coarse moduli space $X$ is an algebraic space of finite type over $k$.
The remainder of this section is devoted to giving a criterion for the properness of $X$ over $k$,
which is described in terms of $\XX$ without making
reference to $X$.
We begin by considering the condition which assures that $X$ is locally separated.

\vspace{2mm}

Let $p$ be a closed point on $\XX$.
Since $p$ is a GIT-like p-stable point,
there exists an effective versal deformation $\xi:\Spec A\to \XX$
such that if $I$ denotes the ideal generated by nilpotent elements of $A$, then there exists a normal subgroup scheme $\mathcal{F}_{\xi}$
of $\Aut_{\XX,\Spec A/I}(\xi|_{\Spec A/I})\to \Spec A/I$ satisfying (b) in Definition~\ref{defgit}.
Let $F$ be a functor on the category of $\Spec A/I\times_k\Spec A/I$-schemes
which to any $(f,g):T\to \Spec A/I\times_k\Spec A/I$
associates the set of sections $\Isom_{\XX,T}(f^*\xi,g^*\xi)/\mathcal{F}_{\xi}(T)$.
Let $\XX_0$ be the reduced stack associated to $\XX$.
Let $\XX_0\to \XX_0^{\textup{rig}}$ be the rigidification associated to an algebraization of $\mathcal{F}_{\xi}$ after replacing $\XX$ by a neighborhood of $p$
(see Proposition~\ref{reducedgit}).
Let $p'$ be the image of $p$ in $\XX\to\XX^{\textup{rig}}$.
Then the composite morphism $\Spec A/I\to \XX_0\to \XX_0^{\textup{rig}}$ is an effective
versal deformation for $p'$ since $\XX_0\to \XX_0^{\textup{rig}}$ is smooth.
By the construction of the rigidification (cf. Remark~\ref{rigrem}),
the functor $F$ is represented by the natural morphism
\[
\Spec A/I\times_{\XX_0^{\textup{rig}}}\Spec A/I\to \Spec A/I\times_{k}\Spec A/I.
\]

Now we prove:

\begin{Proposition}
\label{locallysep}
Suppose that the functor $F$ is proper over $\Spec A/I\times_k\Spec A/I$.
Then $\XX^{\textup{rig}}$ is separated in a neighborhood of $p'$.
\end{Proposition}

\Proof
Without loss of generality, we may assume that $\XX$ is reduced.
Namely, $\XX_0=\XX$ and $I=0$.
It suffices to prove that there exists a smooth neighborhood $U\to \XX^{\textup{rig}}$ of $p'$, such that $U$ is an affine scheme and $U\times_{\XX^{\textup{rig}}}U\to U\times_k U$ is proper.
To show this, we need the following Lemma.

\begin{Lemma}
\label{spreadproper}
Let $X$ and $Y$ be schemes locally finite presentation
over a quasi-compact excellent scheme $S$.
Let $X\to Y$ be a morphism over $S$.
Let $s$ be a closed point,
and let $R:=\hat{\OO}_{S,s}$ be the completion of the local ring at $s$.
Suppose that the morphism
$X\times_{S}\Spec R\to Y\times_{S}\Spec R$ induced by
$\Spec R\to S$ is proper. Then there exists a neighborhood
$W\subset S$ of $s$ such that the induced morphism $X\times_{S}W\to Y\times_SW$
is proper.
\end{Lemma}

\Proof
Let $P$ be a functor on the category of $S$-schemes
which to $T\to S$ associates the set consisting of one element
if $X\times_{S}T\to Y\times_{S}T$ is proper, and associates
the empty set if otherwise.
By Theorem~\ref{limit4}, we see that
this functor is locally of finite presentation.
Then applying Artin's approximation (\cite{A}) we conclude that
there exists an \'etale neighborhood $U\to S$ of $s$,
such that $X\times_{S}U\to Y\times_{S}U$ is proper.
Let $W$ be the image of $U$ in $S$. It is an open set.
By the descent theory, $X\times_{S}W\to Y\times_{S}W$ is proper.
\QED

We continue the proof of Proposition~\ref{locallysep}.
By Lemma~\ref{algebraization},
there exist an affine scheme $U$, a smooth morphism
$w:U\to \XX^{\textup{rig}}$, and a closed point $u\in U$ such that
$A\cong \hat{\OO}_{U,u}$ and $\xi|_{\Spec A}\cong w|_{\Spec A}$
Consider the natural morphism
$U\times_{\XX^{\textup{rig}}}\Spec A\to U\times_k\Spec A$ over $U$.
Then applying Lemma~\ref{spreadproper} to this diagram, we see that
there exists a neighborhood $V$ of $u$ such that
the restriction $V\times_{\XX^{\textup{rig}}}\Spec A\to V\times_k\Spec A$
is proper. (The algebraic space $U\times_{\XX^{\textup{rig}}}\Spec A$
is a scheme since it is quasi-finite and separated over $U\times_k\Spec A$
(cf. \cite[(A.2)]{LM}).
Applying Lemma~\ref{spreadproper} to the morphism
$V\times_{\XX^{\textup{rig}}}U\to V\times_kU$ over $U$ again, we conclude that
after shrinking the neighborhood $V$,
 the morphism $V\times_{\XX^{\textup{rig}}}V\to V\times_kV$ is proper.
 \QED

\begin{Corollary}
\label{locallysep2}
Let $q$ be the image of $p$ on the coarse moduli space.
Under the assumption of Proposition~\ref{locallysep},
$X$ is separated in a neighborhood of $q$.
\end{Corollary}

\Proof
It follows from Keel-Mori theorem since by Proposition~\ref{locallysep} $\XX^{\textup{rig}}$ has finite
diagonal after shrinking $\XX^{\textup{rig}}$.
\QED

We say that $p$ satisfies the locally separated property if
there exist an effective versal deformation $\xi:\Spec A\to \XX$ for $p$
and a closed subgroup scheme $\mathcal{F}_{\xi}$ as above, such that
$F$ (see Proposition~\ref{locallysep}) is proper over $\Spec A/I\times_k \Spec A/I$, where $I$ is the ideal generated by nilpotent elements of $A$.
From now on, we suppose that every point satisfies the locally separated property. Namely, according to Corollary~\ref{locallysep2},
the coarse moduli space $X$ is {\it locally separated}.
Next we consider a criterion for $X$ to be universally closed over $k$.

\begin{Proposition}
\label{univclosed}
The coarse moduli space $X$ is universally closed over $k$ if and 
only if $\XX$ is universally closed over $k$, i.e., $\XX$ satisfies a valuative criterion
in \cite[Th\'eor\`em 7.3]{LM}.
\end{Proposition}

\Proof
The ``only if'' direction is clear since $\pi$ is a universally closed map.
The ``if" direction  follows from the easy fact:
If $\XX$ is universally closed over $k$, and $\XX\to X$ is surjective,
then $X\to \Spec k$ is a universally closed map.
\QED

Finally, we consider a valuative criterion for the separatedness of $X$.

\begin{Lemma}
\label{closeddiagonal}
The coarse moduli space $X$ is separated if and only
if the image of diagonal map $\XX\to \XX\times_k\XX$ is closed.

\end{Lemma}

\Proof
The diagonal $X\to X\times_kX$ is a quasi-compact immersion
since $X$ is locally separated.
Thus, $X$ is separated if and only if the image of $X\to X\times_kX$ is
closed. Suppose that
the image $\mathcal{Z}$ of diagonal map $\XX\to \XX\times_k\XX$ is closed.
Then the the image of $\mathcal{Z}$ in $X\times_kX$ is closed.
Indeed, the composite map $\XX\times_k\XX\to \XX\times_kX\to X\times_kX$
is closed since $\XX\to X$ is universally closed.
Conversely, assume that $X$ is separated.
Let $Z$ be the image of diagonal $X\to X\times_kX$.
Since $\XX\to X$ is a coarse moduli map, thus
the image of diagonal map $\XX\to \XX\times_k\XX$
set-theoretically coincides with the preimage of $Z$ under $\XX\times_k\XX\to X\times_kX$.
It follows that
the image of diagonal map $\XX\to \XX\times_k\XX$
is closed.
\QED

\begin{Proposition}
\label{sepcr}
The following conditions are
equivalent:
\begin{enumerate}
\renewcommand{\labelenumi}{(\roman{enumi})}
\item $X$ is separated over k.
\item Let $R$ be a valuation ring with quotient field
$K$, and let $\alpha,\beta$ be objects in $\XX(R)$ such that
there is an isomorphism $\alpha|_K\cong \beta|_K$.
Then the fiber of $\textup{Isom}_{R}(\alpha,\beta)$ over the closed point
of $R$ is nonempty.
%
%

\end{enumerate}
\end{Proposition}

\Proof
We first show that (ii) implies (i).
Since $X$ is locally separated, thus by Lemma~\ref{closeddiagonal}
it is enough
to prove that the image $\mathcal{Z}$ of the diagonal map
$\Delta:\XX\to \XX\times_k\XX$ is closed.
By \cite[(5.9.4)]{LM}, the underlying set of
$\mathcal{Z}$ is a constructible set.
It suffices to show that $\mathcal{Z}$ is stable under specialization.
Thus it is enough to prove that if $v$ is a generic point on $\XX$,
and $y'$ is a specialization of $y=\Delta(v)$, then there
exists a point $v'\in \XX$ lying over $y'$.
According to \cite[(7.2)]{LM},
there exist
a valuation ring $R$ with quotient field $K$,
$\Spec K\to \XX$ lying over $v$, and 
$\Spec R\to \XX\times_k\XX$
such that the diagram
\[
\xymatrix{
\Spec K \ar[d] \ar[r] & \XX \ar[d]^{\Delta} \\
\Spec R \ar[r] & \XX\times_k\XX \\
}
\]
commutes, and
the closed point of $\Spec R$ maps to $y'$
under $\Spec R\to \XX\times_k\XX$.
Then applying the condition (ii) to $\textup{Isom}_R(\alpha,\beta)\cong \Spec R\times_{(\XX\times_k\XX)}\XX\to \Spec R$
we see that $y'$ belongs to $\mathcal{Z}$.

We next show that (i) implies (ii). 
The image of $(\alpha,\beta):\Spec R\to\XX\times_k\XX$
set-theoretically belongs to $\mathcal{Z}$ since $\mathcal{Z}$ is closed
by Lemma~\ref{closeddiagonal}. 
Hence there exists a point $s\in \textup{Isom}_{R}(\alpha,\beta)$
lying over the closed point of $\Spec R$. \QED

\begin{Remark}
For stable algebraic stacks,
it seems that three
conditions: the locally separated condition, ``universally closed condition''
 (cf. Proposition~\ref{univclosed}) and (ii) in
 Proposition~\ref{sepcr} are important. I like to
 think of three conditions
as ``{\it virtual properness}'' of $\XX$. In other words,
under these conditions, $\XX$ behaves like proper in some contexts.
For example, it is hopeful that the class of such properness provides a
good setting in which one has a good theory of Riemann-Roch.
\end{Remark}

\renewcommand{\thesection}{}
\renewcommand{\theTheorem}{A.\arabic{Theorem}}
\renewcommand{\theClaim}{A.\arabic{Theorem}.\arabic{Claim}}
\renewcommand{\theequation}{A.\arabic{Theorem}.\arabic{Claim}}
\setcounter{Theorem}{0}
\section*{Appendix}

{\it Limit argument}.
Let $S_0$ be a scheme (resp. quasi-separated algebraic space).
Let $I$ be a filtered inductive system in the sense of \cite[section 1]{Art0}
or \cite[Appendix A]{Milne}.
Let us consider functor
\[
I\to (\textup{quasi-coherent}\ \OO_{S_0}\textup{-algebras})
\]
sending
$\alpha$ to $\mathcal{A}_{\alpha}$.
Let \[
S=\displaystyle\lim_{\longleftarrow \atop \alpha\in I}S_{\alpha}
\]
be the associated projective system of schemes (resp. quasi-separated algebraic spaces)
that are affine over $S_0$. (By \cite[Appecdix A, Corollary 2]{Milne}
it is represented by a scheme (resp. quasi-separated algebraic space).)
We would like to recall the some results of limits arguments in \cite[IV (8.2)]{EGA1}.
For our purpose, we need some assertions in the case of algebraic spaces,
though the generalizations are straightforward.
However, we could not find the appropriate literature,
thus we decided that it is best to collect them here.


\begin{Theorem}[cf. EGA IV (8.6.3)]
\label{limitclosed}
Assume that $S_0$ is a quasi-compact and quasi-separated algebraic space.
Let $Y\subset S$ be a closed subspace.
Then there exist $\lambda\in I$ and a closed subspace
$Y_{\lambda}\subset S_{\lambda}$ such that
$S\times_{S_{\lambda}}Y_{\lambda}=Y$.
Moreover, $Y_{\lambda}\subset S_{\lambda}$ is unique in the following sense:
If there exists another closed subscheme $Y_{\lambda'}\subset S_{\lambda'}$
such that $S\times_{S_{\lambda'}}Y_{\lambda'}=Y$,
then there exist $\lambda\to \mu$ and $\lambda'\to \mu$
such that $Y_{\lambda}\times_{S_{\lambda}}{S_{\mu}}=Y_{\lambda'}\times_{S_{\lambda'}}S_{\mu}$.
\end{Theorem}

\Proof
Assume $S_0$ is an quasi-compact and quasi-separated algebraic space.
Let $\tilde{S}_0\to S_0$ be an \'etale surjective morphism from a quasi-compact and separated scheme $\tilde{S}_0$.
Considering the base changes $\{\tilde{S}_0\times_{S_0}S_{\alpha}\}_{\alpha\in I}$
the uniqueness follows from the case when $S_0$ is a scheme.
Next we will prove the existence.
Let $X_{\alpha}:=\tilde{S}_0\times_{S_0}S_{\alpha}$
and $X=\tilde{S}_0\times_{S_0}S$. Let $p_{\alpha}:X_{\alpha}\to S_{\alpha}$
and $p:X\to S$ be natural projections.
Then there exist $\lambda\in I$ and a closed subscheme $W_{\lambda}\subset X_{\lambda}$ such that $X\times_{X_{\lambda}}W_{\lambda}=p^{-1}(Y)$.
Let $R_{\alpha}:=X_{\alpha}\times_{S_{\alpha}}X_{\alpha}$
and $R:=X\times_{S}X$.
Since $X_{\alpha}$ is a quasi-compact and quasi-separated scheme, thus $R_{\alpha}$ is also a quasi-compact and quasi-separated scheme.
Let $\textup{pr}_1,\textup{pr}_2:R_{\alpha}\rightrightarrows X_{\alpha}$
be the first and second projection respectively. (Here we abuse notation and omit the index ``$\alpha$''.)
It suffices to show that there exists some arrow $\lambda\to \mu$
such that $\textup{pr}_1^{-1}(W_{\mu})=\textup{pr}_2^{-1}(W_{\mu})$,
where $W_{\mu}=W_{\lambda}\times_{X_{\lambda}}X_{\mu}$.
Since $\textup{pr}_1^{-1}(p^{-1}(Y))=\textup{pr}_2^{-1}(p^{-1}(Y))$,
the uniqueness $\textup{pr}_1^{-1}(W_{\lambda})$ and $\textup{pr}_2^{-1}(W_{\lambda})$ in the system $\{R_{\alpha}\}_{\alpha\in I}$ implies that
there exists $\lambda\to \mu$ such that $\textup{pr}_1^{-1}(W_{\mu})=\textup{pr}_2^{-1}(W_{\mu})$. Thus $W_{\mu}$ descends to a closed subspace $Y_{\mu}\subset S_{\mu}$ such that $S\times_{S_{\mu}}Y_{\mu}=Y$.
\QED

\begin{Theorem}[cf. EGA IV (8.10.5) (17,7.8)]
\label{limit4}
Let $S_0$ be a quasi-compact scheme.
Let $X_{\lambda}$ and $Y_{\lambda}$ be schemes of finite presentation
over $S_{\lambda}$.
Let $X_{\lambda}\to Y_{\lambda}$ be an $S_{\lambda}$-morphism.
Consider the following properties: (i) an isomorphism, (ii) a closed immersion,
(iii) quasi-finite,
(iv) finite, (v) affine, (vi) flat, (vii) smooth, (viii) proper, (ix) separated.
If $X_{\lambda}\times_{S_{\lambda}}S\to
Y_{\lambda}\times_{S_{\lambda}}S$ has one or more properties (i)--(viii), then there exists
an arrow $\lambda\to \mu$ such that $X_{\lambda}\times_{S_{\lambda}}S_{\mu}\to
Y_{\lambda}\times_{S_{\lambda}}S_{\mu}$ has the same properties.
Furthermore the similar assertion of (ix) for an algebraic space $X_{\lambda}$ holds.

\end{Theorem}

\Proof
What we have to prove is the last statement.
It suffices to show that
$X_{\mu}\to X_{\mu}\times_{Y_{\mu}}X_{\mu}$
is a closed immersion for some $\mu$.
Note that $X_{\lambda}\to X_{\lambda}\times_{Y_{\lambda}}X_{\lambda}$
is quasi-affine. Thus our claim follows from (ii) by descent theory.
\QED

\begin{Proposition}
\label{limit2}
Suppose that $S_{\alpha}$ is a noetherian scheme for any $\alpha\in I$.
Let $f:X_{\lambda}\to S_{\lambda}$ be an algebraic space of finite type
over $S_{\lambda}$. Consider the following properties:
(i) affine,
(ii) flat,
(iii) finite,
(iv) smooth.
If $X:=X_{\lambda}\times_{S_{\lambda}}S\to
S$ has one or more properties (i)--(iv), then there exists
an arrow $\lambda\to \mu$ such that $X_{\mu}:=X_{\lambda}\times_{S_{\lambda}}S_{\mu}\to
S_{\mu}$ has the same properties.
\end{Proposition}

\Proof
We first prove (ii).
Let $W_{\alpha}\to X_{\alpha}$ be an \'etale surjective morphism from a quasi-compact
scheme $W_{\alpha}$. Note that $W_{\alpha}\to X_{\alpha}$ is a quasi-compact morphism
since $X_{\alpha}$ is quasi-separated.
By the descent theory, it is enough to show that
there exists an arrow $\alpha\to \mu$ such that
$W_{\alpha}\times_{S_{\alpha}}S_{\mu}$ is flat over $S_{\mu}$.
Applying Theorem~\ref{limit4} to the projective system $\{W_{\alpha}\times_{S_{\alpha}}S_{\lambda}\}_{\lambda\in I/\alpha}$,
we obtain (i). (Here $I/\alpha$ is the inductive system over $\alpha$.)

Next we prove (iv).
Let $[R_{\alpha}=W_{\alpha}\times_{X_{\alpha}}W_{\alpha}\rightrightarrows W_{\alpha}]$ denote the \'etale equivalence relation for $X_{\alpha}$.
Since $X\to S$ is finite, in particular, separated, thus
by Theorem~\ref{limit4} there exists an arrow $\alpha\to \mu$ such that
$(\textup{pr}_1,\textup{pr}_2):R_{\mu}\to W_{\mu}\times_{S_{\mu}} W_{\mu}$ is a closed immersion,
that is, $X_{\mu}$ is separated over $S_{\mu}$. Furthermore, again by Theorem~\ref{limit4} we may assume that $W_{\mu}$ is quasi-finite over $S_{\mu}$,
that is, $X_{\mu}$ is quasi-finite over $S_{\mu}$.
Then by \cite[(A.1)]{LM}, $X_{\mu}$ is a scheme.
Therefore we have (ii) by Theorem~\ref{limit4}.

Next we prove (iii).
According to Theorem~\ref{limit4} (iii) we may assume that
$f$ is quasi-finite.
There exists some $\mu$ such that $f_{\mu}$ is separated by Theorem~\ref{limit4}.
Then by ****** $X_{\mu}$ is a scheme. Therefore the scheme case implies
(iii).

Finally, we prove (i).
Since $X=X_{\lambda}\times_{S_{\lambda}}S$ is affine of finite type over $S$,
there exists a closed immersion $\iota:X\to \Spec S[t_1,\ldots,t_n]$.
It is enough to show that there is a closed immersion
$\iota_{\mu}:X_{\mu}\to \Spec S_{\mu}[t_1,\ldots,t_n]$ which
induces $\iota$.
Let $[R_{\lambda}=W_{\lambda}\times_{X_{\lambda}}W_{\lambda}\rightrightarrows W_{\lambda}]$ denote an \'etale equivalence relation for $X_{\lambda}$.
Using Proposition~\ref{homlimit} and the
\'etale equivalence relation, there exists
$X_{\mu}\to \Spec S_{\mu}[t_1,\ldots,t_n]$ for some $\mu$.
Then by the above (iii) and Theorem~\ref{limit4} (ix),
we may suppose that $X_{\mu}\to \Spec S_{\mu}[t_1,\ldots,t_n]$
is quasi-finite and separated. In particular, we can assume
$X_{\mu}$ is a scheme. Now our claim follows from Theorem~\ref{limit4}
(ii).
\QED

\begin{Proposition}
\label{limitconnected}
Suppose that $S_{\alpha}$ is a noetherian scheme for any $\alpha\in I$.
Let $f_{\lambda}:X_{\lambda}\to S_{\lambda}$ is a scheme of finite type over $S_{\lambda}$.
Suppose that all geometric fibers of $X_{\lambda}\times_{S_{\lambda}}S\to S$
are connected.
Then there exists an arrow $\lambda\to \mu$
such that all geometric fibers of $X_{\lambda}\times_{S_{\lambda}}S_{\mu}\to S_{\mu}$
are connected.

\end{Proposition}

\Proof
Let $E$ be the set consisting of points $s\in S_{\alpha}$,
such that $f^{-1}(s)$ is geometrically connected.
By \cite[IV, 9.7.7]{EGA1}, the set $E$ is a
constructible set in $S_{\lambda}$.
Moreover, if a point $s\in S_{\lambda}$ does not lie in the image of $S\to S_{\lambda}$,
then there exists $\lambda\to \lambda'$ such that $s$ does not lie in
the image of
$S_{\lambda'}\to S_{\lambda}$. On the other hand, the image of $S\to S_{\lambda}$
is contained in $E$.
Therefore, using noetherian induction argument, we easily see that
there exists $\lambda \to\mu$ such that the image of $S_{\mu}\to S_{\lambda}$
is contained in $E$. This completes the proof.
\QED

\begin{Theorem}[cf. EGA IV (8.8.2)]
\label{homlimit}
Assume that $S_0$ be a quasi-compact and quasi-separated scheme.
Let $X_{\alpha}$ be a quasi-compact scheme and let
$Y_{\alpha}$ be a scheme of finite presentation over $S_{\alpha}$.
Let $X_{\mu}=X_{\alpha}\times_{S_\alpha}S_{\mu}$ and $Y_{\mu}=Y_{\alpha}\times_{S_\alpha}S_{\mu}$.
Let $X=X_{\alpha}\times_{S_\alpha}S$ and $Y=Y_{\alpha}\times_{S_\alpha}S$.
Then the natural map
\[
A=\displaystyle\lim_{\longrightarrow \atop \lambda} \Hom_{S_{\lambda}}(X_{\lambda},Y_{\lambda})\to \Hom_{S}(X,Y)
\]
is bijective. 
\end{Theorem}

\begin{Proposition}
\label{grouplimit}
Assume that $S_0$ be a quasi-compact and quasi-separated scheme.
Let $G_{\alpha}\to S_{\alpha}$ be a separated group algebraic space
of finite presentation over $S_{\alpha}$.
Let $H_{\alpha}\subset G_{\alpha}$ be a closed subspace
such that $H=H_{\alpha}\times_{S_{\alpha}}S$ is a subgroup space of $G=G_{\alpha}\times_{S_{\alpha}}S$.
Then there exists an arrow $\alpha\to \mu$ such that
$H_{\mu}:=H_{\alpha}\times_{S_{\alpha}}S_{\mu}$ is a subgroup space of $G_{\mu}=G_{\alpha}\times_{S_{\alpha}}S_{\mu}$.
If $H$ is normal in $G$,
then $H_{\alpha}\times_{S_{\alpha}}S_{\mu}$ can be chosen to be normal.
\end{Proposition}

\Proof
It suffices to show that there exists an arrow $\alpha\to \mu$
such that $H_{\mu}$ has the following properties:
(i) the multiplication $m:G_{\mu}\times_{S_{\mu}}G_{\mu}\to G_{\mu}$
induces $H_{\mu}\times_{S_{\mu}}H_{\mu}\to H_{\mu}$, (ii) the inverse
$i:G_{\mu}\to G_{\mu}$ carries $H_{\mu}$ to $H_{\mu}$,
(iii) the unit section $e:S_{\mu}\to G_{\mu}$ factors through
$H_{\mu}$.
Let us consider the property (i).
Let 
\[
\xymatrix{
U_{\alpha} \ar[r]^{\tilde{m}} \ar[d]^{p} & V_{\alpha} \ar[d]^{q}\\
G_{\alpha}\times_{S_{\alpha}}G_{\alpha} \ar[r]^(0.7){m} & G_{\alpha} \\
}
\]
be a commutative diagram where $U_{\alpha}$ and $V_{\alpha}$
are
quasi-compact schemes and two vertical arrows are \'etale surjective morphisms.
 (Such a diagram exists.)
It suffices to prove that after the base change by some arrow $\alpha\to \mu$
the pullback $p^{-1}(H_{\mu}\times_{S_{\mu}}H_{\mu}$) maps to
$q^{-1}(H_{\mu})$.
Since $H$ is a subgroup of $G$, in  particular
$G\times_{S}G\to G$ induces $H\times_{S}H\to H$, thus applying Theorem~\ref{homlimit}
to $p^{-1}(H_{\alpha}\times_{S_{\alpha}}H_{\alpha})\to V_{\alpha} \leftarrow q^{-1}(H_{\alpha})$ we see that there exists an arrow $\alpha\to \mu$
such that
$p^{-1}(H_{\mu}\times_{S_{\mu}}H_{\mu})\to V_{\mu}$ factors through
$q^{-1}(H_{\mu})\subset G_{\mu}$.
By a similar argument, we may assume that $H_{\mu}$ has also the property (ii).

Next, consider the property (iii). Let
\[
\xymatrix{
W_{\alpha} \ar[r]^{\tilde{e}} \ar[d]^{r} & V_{\alpha} \ar[d]^{q}\\
S_{\alpha}\ar[r]^{e} & G_{\alpha} \\
}
\]
be a commutative diagram where $W_{\alpha}$ and $V_{\alpha}$
are
quasi-compact schemes and two vertical arrows are \'etale surjective morphisms.
It suffices to prove that after the base change by some arrow $\alpha\to \mu$
the unit section $W_{\mu}\to V_{\mu}$ factors through $q^{-1}(H_{\mu})$.
Again by applying Theorem~\ref{homlimit} to $W_{\alpha}\to V_{\alpha}\leftarrow
q^{-1}(H_{\alpha})$, we conclude that there exists an arrow $\alpha\to \mu$
such that
$W_{\mu}\to V_{\mu}$ factors through
$q^{-1}(H_{\mu})\subset G_{\mu}$.
By what we have proven,
we conclude that
there exists an arrow $\alpha\to \mu$ such that
$H_{\mu}$ is a subgroup space of $G_{\mu}$.

Finally, we will prove the last assertion.
Consider the morphism
\[
\phi:G_{\alpha}\times_{S_{\alpha}}G_{\alpha}\to G_{\alpha}
\]
which sends $(g_1, g_2)$ to $g_1g_2g_1^{-1}$.
A group subspace $H_{\alpha}$ is normal in $G_{\alpha}$
if and only if $\phi(G_{\alpha}\times_{S_{\alpha}}H_{\alpha}) \subset H_{\alpha}$ (scheme-theoretically).
Thus by a similar argument using Theorem~\ref{homlimit} to the above we easily see the last assertion.
\QED

\vspace{5mm}

{\it Rigidification}.
For the convenience, we would like to recall the rigidifications of algebraic stacks, that has been
discussed by many authors (see for example \cite{LM}, \cite{AGV}, \cite{AOV},\cite{Ro}).

Let $\XX$ be a quasi-separated algebraic stack locally of finite presentation over a quasi-separated scheme $S$.
Let $\GGG$ be a closed subgroup stack in the inertia stack $I\XX\to \XX$,
that is flat and of finite presentation over $\XX$.
Namely, the multiplication $I\XX\times_{\XX}I\XX\to I\XX$ induces $\GGG\times_{\XX}\GGG\to\GGG$,
the inverse $I\XX\to I\XX$ sends $\GGG$ to $\GGG$,
and in addition the unit section $\XX\to I\XX$ factors through $\GGG$.
Let $V$ be an affine $S$-scheme. Let $\xi\in\XX(V)$
be an object
and $V\to \XX$ the corresponding morphism.
It gives rise to a natural morphism $h:\Aut_{\XX,V}(\xi)\to I\XX$.
The inverse image $\GGG(\xi):=h^{-1}(\GGG)\subset \Aut_{\XX,V}(\xi)$ is a normal closed
subgroup. (Every automorphism $\sigma:\xi\to\xi$ gives rise to
an inertia automorphism $\Aut_{\XX,V}(\xi)\to\Aut_{\XX,V}(\xi)$ and
it sends $h^{-1}(\GGG)$ to $h^{-1}(\GGG)$. Thus $\GGG(\xi)$ is normal.)

\begin{Theorem}
\label{rigidification}
There exist an algebraic stack $\YY$ locally of finite
presentation and a morphism $f:\XX\to \YY$ such that:
\begin{enumerate}
\renewcommand{\labelenumi}{(\roman{enumi})}

\item $f:\XX\to\YY$ is an fppf gerbe.

\item For any affine $S$-scheme $V$ and any object $\xi\in \XX(V)$,
the homomorphism of group algebraic spaces
\[
\Aut_{\XX,V} (\xi)\to \Aut_{\YY,V} (f(\xi))
\]
is surjective and its kernel is $\GGG(\xi)$.

\item Let $V\to \YY$ be a morphism from an affine scheme $V$
and let $\xi\in\XX(V)$ be the corresponding object.
Then there exists a natural isomorphism $V\times_{\YY}\XX\cong B_{V}\GGG(\xi)$
over $V$, where $B_V\GGG(\xi)$ is the classifying stack of $\mathcal{G}(\xi)$.

\item Let $g:\XX\to \WW$ be a morphism of algebraic stacks 
such that for any object $\xi\in\XX$, $\GGG(\xi)$ lies in the kernel of
$\Aut_{\XX,V}(\xi)\to \Aut_{\WW,V}(g(\xi))$.
Then there exists a morphism $h:\YY\to \WW$ such that
$h\circ f$ is isomorphic to $g$. It is unique up to isomorphism.

\item $\YY$ admits a coarse moduli space if and only if $\XX$ has one.
If they have, both coarse moduli spaces coincide.

\item If $\XX$ is a Deligne-Mumford stack, then $\YY$ is so.

\end{enumerate} 

\end{Theorem}

\Proof
Assertions (i), (ii) (iii) are proved in \cite[(A.1)]{AOV}.
To see (iv), (v) and (vi), we should go back to the construction
of $\YY$.
Let $\YY^p$ be the fibered category defined as follows.
The objects of $\YY$ are those of $\XX$.
Let $\phi:V'\to V$ be a morphism of affine $S$-schemes.
For any $\xi'\in\XX(V')$ and any $\xi\in\XX(V)$,
we define the set of morphisms $\Hom_{\YY^p}(\xi',\xi)$ over $V'\to V$
to be the set of sections of the quotient
\[
\underline{\operatorname{Isom}}_{\XX,V'}(\xi',\phi^*\xi)\slash\GGG(\xi')
\]
where the action of $\GGG(\xi')$ is the natural faithful right action
defined by the composition.
Then it is not hard to show that $\YY^p$ is a prestack.
Let $\YY$ be a stack associated to the prestack $\YY^p$.
It is showed in \cite[(A.1)]{AOV} that $\YY$ is an algebraic stack locally
of finite presentation over $S$ with properties (i), (ii) and (iii).
By the construction of $\YY^p$, clearly,
the natural morphism $\XX\to\YY^p$ is universal for maps from $\XX$
to prestacks, that kill $\GGG(\xi)$ for any object $\xi\in \XX$.
Since $\YY$ is the associated stack, thus (iv) follows.
The ``if part'' of (v) follows from (iv) and (i).
(Note that for any algebraically closed $S$-field $K$,
$\XX\to\YY$ induces a bijective map $\XX(K)\slash\sim\to\YY(K)/\sim$.
Here $\sim$ means ``up to isomorphism''.)
The ``only if'' part also follows from (iv).
The last assertion in (v) is clear.
To see (vi), notice that for each object $\xi\in\XX(V)$
the group space $\Aut_{\XX,V}(\xi)\to V$ is unramified because $\XX$ is a Deligne-Mumford stack (cf. \cite[(4.2)]{LM}).
For each object $\xi\in\XX(V)$, $\GGG(\xi)$ is flat and thus $\GGG(\xi)$ is \'etale over $V$.
By (iii), we see that $\XX\to \YY$ is \'etale. Hence $\YY$ is a Deligne-Mumford stack.
\QED

\begin{Remark}
\label{rigrem}
The morphism $\XX\to \YY$ is characterized by the universal 
property in (iv). We refer to this morphism as the rigidification of $\XX$
(or rigidifying morphism)
with
respect to $\GGG$.
If objects $\xi,\eta \in\YY(V)$ arise from $\XX(V)$,
then the above proof reveals that there exists an isomorphism
\[
\Isom_{\YY,V}(\xi,\eta)\cong \Isom_{\XX,V}(\xi,\eta)/\GGG(\xi). 
\]
Also, the algebraic space on the right side is isomorphic to
$\GGG(\eta)\setminus\Isom_{\XX,V}(\xi,\eta)$.

\end{Remark}

\end{document}